\magnification=\magstephalf

\input amstex
\input epsf
\documentstyle{amsppt}

\hsize=16truecm
\vsize=21.5truecm

\def\({\left(\,}
\def\){\,\right)}
\def\[{\left[\,}
\def\]{\,\right]}

\def\<{\left\langle}
\def\>{\right\rangle}

\def\Var{\text{Var}}


\def \eqd{\buildrel d \over =}
\def \convd{\buildrel d \over \longrightarrow}

\def\a{\alpha}
\def\b{\beta}
\def\d{\delta}
\def\e{\varepsilon}
\def\g{\gamma}

\def\l{\lambda}

\def\z{\zeta}
\def\cd{\Delta}

\def\bZ{{\bold Z}}

\def\A{{\Cal A}}

\def \Ga{\Gamma}


\def\sp{\vskip1ex}
\def\z{\zeta}
\def\ps{\psi}
\def\inv{^{-1}}
\def\ra{\rightarrow}
\def\l{\ell}

\def\iy{\infty}
\def\be{\begin{equation}}
\def\ee{\end{equation}}
\def\ov{\over}
\def\al{\alpha}
\def\noi{\noindent}
\def\cd{\cdots}
\def\dl{\delta}
\def\ph{\varphi}
\def\Pr{P}

\def\si{\sigma}

\def\ub{\bar{u}} 
\def\lan{\left\langle} \def\ran{\right\rangle}



\hyphenation{per-co-la-tion}
\hyphenation{per-co-la-ting}
\hyphenation{Wads-worth}
\hyphenation{Green-berg}
\hyphenation{Has-tings}
\hyphenation{pub-li-cation}

\newdimen\howmuch

\howmuch=0pt

\def \centerbox#1{\setbox0=#1
\advance\howmuch by \the\wd0
\multiply\howmuch by -1
\advance\howmuch by \hsize
\divide\howmuch by 2
\moveright\howmuch #1}

\def\beginitems{\begingroup
    \medbreak
    \parskip=5pt
    \advance \parindent by 2em
    \vskip-\parskip}
\def\enditems{
    \medbreak
    \vskip-\parskip
    \endgroup}

\font\sc=cmcsc10


\document


\pageno=0

\footline={\hfil}

\null
\vskip0.5cm
 
\centerline{\bf FLUCTUATIONS IN THE COMPOSITE REGIME}
\centerline{\bf OF A DISORDERED GROWTH MODEL}  

\vskip1cm  

\centerline{\sc Janko Gravner}
\centerline{\rm Department of Mathematics }
\centerline{\rm University of California}
\centerline{\rm Davis, CA 95616}
\centerline{\rm email: \tt gravner\@math.ucdavis.edu}

\vskip 0.2cm

\centerline{\sc Craig A. Tracy}
\centerline{\rm Department of Mathematics}
\centerline{\rm Institute of Theoretical Dynamics}
\centerline{\rm University of California}
\centerline{\rm Davis, CA 95616}
\centerline{\rm email: \tt tracy\@itd.ucdavis.edu}

\vskip 0.2cm

\centerline{\sc Harold Widom}
\centerline{\rm Department of Mathematics}
\centerline{\rm University of California}
\centerline{\rm Santa Cruz, CA 95064}
\centerline{\rm email: \tt widom\@math.ucsc.edu}

\vskip 0.2cm

\centerline{(Version 2, March 22, 2002)}

\vskip 0.2cm

\flushpar{\bf Abstract.} We continue to study a model of 
disordered interface 
growth in two dimensions. The 
interface is given by a height function 
on the sites of the one--dimensional integer lattice and grows  
in discrete time: 
(1) the height 
above the site $x$ adopts the height above the 
site to its left 
if the latter height is larger, (2) otherwise, the height above $x$ 
increases by 1 with probability $p_x$. We assume that 
$p_x$ are chosen independently at random with a common 
distribution $F$, and that the initial state is such that 
the origin is far above the other sites. Provided that the tails 
of the distribution $F$ at its right edge are sufficiently 
thin, there exists  a nontrivial composite regime in which the 
fluctuations of this interface 
are governed by extremal statistics of $p_x$. 
In the quenched case, the said fluctuations are 
asymptotically normal, while in the annealed case they satisfy
the appropriate extremal limit law. 
  
\vskip 0.2cm

\flushpar 2000 {\it Mathematics Subject Classification\/}. Primary 60K35.
Secondary 05A16, 33E17, 60K37, 60G70, 82C44.  

\vskip0.2cm

\flushpar {\it Keywords\/}: growth model,  
fluctuations, Fredholm determinant, phase transition, 
saddle point analysis, extremal order statistics.

\vskip0.2cm 

\flushpar {\bf Acknowledgments.} 
This work was partially supported by 
National Science Foundation  grants DMS--9703923, DMS--9802122, and 
DMS--9732687, as well as
the Republic of Slovenia's Ministry of Science Program Group 503.
Special thanks go to Harry Kesten, who supplied the main idea for the proof 
of Lemma~6.1. The authors are also thankful to the referee
for the careful reading of the manuscript and suggestions for its
improvement.


\vfill\eject

\baselineskip=15pt

\parskip=13pt 

\pageno=1

\centerline{\bf FLUCTUATIONS IN THE COMPOSITE REGIME}
\centerline{\bf OF A DISORDERED GROWTH MODEL} 
   
\vskip0.3cm  

\centerline{\sc Janko Gravner,  Craig A. Tracy, Harold Widom}

\vskip0.5cm

\subheading{1. Introduction}

Disordered systems, which are, especially in the 
context of magnetic materials, often referred to as 
{\it spin glasses\/}, have been the subject 
of much research since the pioneering work in the 1970s. 
The vast majority of this work is
nonrigorous, based on simulations and techniques 
for which a proper mathematical foundation is yet 
to be developed. (See [MPV] for early developments and 
[Tal] for a nice overview 
of the mean field approach.)  As a result, there is  
a large number of new and 
intriguing phenomena 
observed in these models which await
rigorous treatment. 
Among the most fundamental of issues are 
the existence and the nature of a phase transition 
into a {\it glassy\/} or  
{\it composite\/} phase: below a critical temperature, the dynamics of a
strongly disordered system becomes extremely slow 
with strong correlations, aging and localization effects 
and possibly many local equilibria. 
We refer the reader to [NSv] and [BCKM] and 
other papers in the same volume for  
reviews and pointers to the voluminous literature and to [NSt1] and [NSt2] 
for some recent rigorous results.
In view of the 
difficulties associated with a detailed 
understanding of realistic spinglass systems, 
other disordered models have been introduced, which 
are more amenable to existing probabilistic methods. 

One of the most successful of such (deceptively) 
simple models is the one--dimensional random walk with random 
rates [FIN1]. In this model, the walker waits at a site $x\in \bZ$ for 
an exponential time with mean $\tau_x$ before jumping to
either of its two neighbors with equal 
probability. The disorder variables $\tau_x$ 
are i.i.d.\ and quenched, that is, 
chosen at the beginning. Provided that the distribution 
of $\tau_x$ has sufficiently fat tails, namely, if 
$P(\tau_x\ge t)$ decays 
for large $t$ as 
as $t^{-\a}$ with $\a<1$, the walk 
exhibits aging and localization effects ([FIN1], [FIN2]). 
Various one--dimensional voter models and stochastic Ising models at 
zero temperature 
can be explicitly represented with random walks. This connection 
has been explored to demonstrate glassy phenomena such as 
aging and chaotic time dependence ([FIN1], [FINS]). 
The positive temperature 
versions of such results remain open problems, 
even in one dimension. 

In contrast with models which are exactly solvable 
in terms of random walks and are by now a classical 
subject in spatial processes ([Gri1], [Lig]), 
techniques based on the RSK algorithm
and random matrix theory have entered into the study 
of growth processes only recently ([BDJ], [Joh1], [Joh2], [BR],
[PS], [GTW1]). 
The purpose of this paper is to employ these 
new methods to 
prove  the existence of a
{\it pure phase\/} and a {\it composite phase\/} in a
disordered growth model. It has been
observed before in similar models 
[SK] that the role of temperature is 
for flat interfaces apparently played by their {\it slope\/}. 
In our case, the initial set is very far from flat 
and ``temperature'' is measured instead by the macroscopic 
direction (from the origin) of points on the boundary. 
We identify precisely the critical direction and 
demonstrate that 
the fluctuations asymptotics provide an order parameter 
that distinguishes  the two phases. 
We emphasize that a hydrodynamic quantity, the asymptotic shape,
has a discontinuity of the first derivative at the transition
point, at which the shape changes from curved to flat. 
However, this  does not signify 
the existence of a new phase as kinks are common in many 
random growth models [GG], 
thus a finer resolution is necessary. 

The particular model we investigate is
{\it Oriented 
Digital Boiling (ODB)\/} (Feb.~12, 1996, Recipe at 
[Gri2], [Gra], [GTW1], [GTW2]), arguably the 
simplest interacting model for a growing interface in the
two--dimensional lattice $\bZ^2$. 
The occupied set, which changes in 
discrete time $t=0,1,2,\dots$, is given 
by $\A_t=\{(x,y): x\in \bZ, y\le h_t(x)\}$. The initial 
state is a long stalk at the origin: 
$$
h_0(x)=\cases
0, &\text{if }x=0,\\
-\infty, &\text{otherwise, }
\endcases
$$
while the time evolution of the height function $h_t$ 
is determined thus: 
$$
h_{t+1}(x)=\max\{h_{t}(x-1), h_{t}(x)+\e_{x,t}\}. 
$$
Here $\e_{x,t}$ are independent Bernoulli random 
variables, with $P(\e_{x,t}=1)=p_x$. Although this model 
is simplistic, note that it does involve 
the roughening noise (random increases) as well 
as the smoothing surface tension effect (neighbor interaction),
the basic 
characteristics of many growth and deposition processes. (See Sections
5.1, 5.2 and 5.4 of [Mea] for an overview of simple models 
of ODB type as well as some other disordered growth processes.)

We will assume, throughout this paper, that the 
disorder variables 
$p_x$ are initially chosen 
at random, independently with a common distribution 
$F(s)=P(p_x\le s)$.   
We use $\<\,\cdot\,\>$ to denote integration with 
respect to $dF$ and label by $p$ a generic random variable 
with distribution $F$. 

It quickly turns out ([GTW1]), that fluctuation in ODB can 
be studied via equivalent increasing path
problems. Start by 
constructing a random $m\times n$ matrix $A=A(F)$, with 
independent Bernoulli entries $\e_{i,j}$ and such 
that $P(\e_{i,j}=1)=p_j$, where, again, $p_j\eqd p$ are i.i.d.
Label columns as usual, 
but rows started at the bottom. We call a sequence 
of 1's in $A$ whose positions have column index nondecreasing and row index 
strictly increasing an {\it increasing path\/} in $A$, and denote  
by $H=H(m,n)$ the length of the longest  
increasing path. Then, 
under a simple coupling, $h_t(x)=H(t-x,x+1)$ ([GTW1]). 
Thus we will
concentrate our attention on the random 
matrix $A$ rather than the associated growth model.
From now on we will also replace $p_i$ with its {\it ordered\/} 
sample, so 
that $p_1\ge p_2\ge \dots \ge p_n$ (see section 2.2 of [GTW1]). 

We initiated the study of ODB in a 
random environment in an earlier paper ([GTW2]), 
from which we now summarize the notation and the 
main results. Throughout, we denote by $b$ 
the right edge of the support of $dF$ and 
assume it is below 1, i.e., 
$$b=\min\{s:F(s)=1\}<1.$$
Moreover, we fix an $\a>0$ and assume that 
$n=\a m$. (Actually, 
$n=\lfloor \a m\rfloor$, but we  
omit the obvious integer parts.) 
As mentioned above, we can expect different behaviors 
for different slopes on the boundary of the 
asymptotic shape, which translates to 
different $\a$'s. To be more precise, we  
define the following critical values 
$$
\aligned
&\a_c=\< \frac{p}{1-p}\>^{-1},\\
&\a_c'=\<\frac{p(1-p)}{(b-p)^2}\>^{-1}.
\endaligned 
$$
Note that the  second critical value is nontrivial, i.e., $\a_c'>0$, 
iff $\<(b-p)^{-2}\><\infty$. 
Next, define 
$c=c(\a,F)$ to be the time constant, 
$\displaystyle
c=c(\a,F)=\lim_{m\to\infty}{H}/m, 
$
which determines the limiting shape 
of $\A_t$, namely $\lim \A_t/t$, as $t\to\infty$.
In Theorem 1 of [GTW2], it was found that $c$ exists 
a.s.\ and is given by 
$$
c(\a,F)=
\cases
b+\a (1-b)\<p/(b-p)\>,&\text{ if }\a\le \a_c', \\
a+\a (1-a)\<p/(a-p)\>,&\text{ if }\a_c'\le \a\le \a_c,\\
1,&\text{ if }\a_c\le \a. 
\endcases
$$
Here $a=a(\a,F)\in [b,1]$  is the unique solution to
$\a\<{p(1-p)}/{(a-p)^2}\>=1.$

In [GTW2], we also determined 
fluctuations in the {\it pure\/} regime $\a_c'<\a< \a_c$. 
(The  {\it deterministic\/} regime $\a_c<\a$ has no 
fluctuations.) The {\it annealed\/} fluctuations ([GTW2], Theorem 2) 
about the 
deterministic shape $c$ grow as $\sqrt m$ and 
are asymptotically normal: 
$$
\frac{H-cm}{\tau_0\sqrt\a \cdot m^{1/2}}\convd N(0,1)
$$
as $m\to\infty$, where 
$
\tau_0^2=\Var({(1-a)p}/{(a-p)}). 
$

By contrast, {\it quenched\/} fluctuations conditioned on the state of the
environment grow more slowly, as $m^{1/3}$, and satisfy 
the $F_2$--distribution known from random matrices ([TW1],
[TW2]). 
To formulate this result, we let $r_j=p_j/(1-p_j)$, define 
$u_n$ to be the solution of 
$${\al\ov n}\,\sum_{j=1}^n{r_j\ov
(1+r_ju)^2}={1\ov(u-1)^2}\tag{1.1}$$
which lies in 
in the interval $(-r_1\inv,\,0)$. This solution  exists provided that $\a 
n^{-1}\sum_{j=1}^n r_j<1$ 
which holds a.s.\ for large $n$ as soon as $\a<\a_c$. 
Next, set $c_n=c(u_n)$ where 
$$c(u)={1\ov 1-u}-{\al\ov n}\sum_{j=1}^n{r_ju\ov 1+r_ju}.\tag{1.2}$$
Then ([GTW2], Theorem 3) 
there exists a constant $g_0\ne 0$ so that 
$$
P\(\frac{H-c_n m}{g_0^{-1} m^{1/3}}\le s\,\mid\, p_1,\dots,p_n\)\to F_2(s), 
$$
as $m\to \infty$, almost surely, for any fixed $s$.

For fluctuation results in this paper we need to impose 
some additional 
assumption on $F$, which are best expressed in terms 
of $G(x)=1-F((b-x)-)$, the distribution function for $b-p$. 
First we list our weaker conditions: 

\beginitems 

\item{(a)} If $x,y\to 0$ and $x\sim y$, then $G(x)\sim G(y)$. 

\item{(b)} If $x,y\to 0$ and $x=O(y)$, then $G(x)= O(G(y))$.

\item{(c)} As $x\to 0$, $G(x)=o(x^2/\log x^{-1})$. 

\enditems

\flushpar Our stronger assumptions on $F$ require that there exists a $\g>0$ 
so that:

\beginitems 

\item{(a$'$)} The 
function $G(x)/x^{\g}$ is nonincreasing in a neighborhood 
of $x=0$. 

\item{(b$'$)} $G(x)=O(x^{2}/\log^\nu x^{-1})$ as $x\to 0$ 
for some $\nu>2\g+4$.  
 
\enditems

If $\a_c'>0$, then automatically $G(x)=o(x^2)$ as $x\to 0$. 
The stronger assumptions thus do not require 
much more: for nicely behaved $G$ they amount to
$G(x)=O(x^2/\log^\nu x^{-1})$ for some 
$\nu>8$. The quenched and annealed fluctuations 
are now determined by the next two theorems. 

\proclaim{Theorem 1} Assume that $0<\a<\a_c'$, let 
$$
\tau^2= {b(1-b)}\(\frac 1\a-\frac 1{\a_c'}\),
$$ 
and let $\Phi$ be the standard normal distribution function. 
If (a)--(c) hold, then for any fixed $s$, as $m\to\infty$,
$$
P\(\frac {H-c_n m+2\tau \sqrt n}{\tau \sqrt n}\le s
\,\mid\, p_1,\dots,p_n\)\to \Phi(s). 
$$
Here, the convergence is in probability if (a)--(c) hold, and 
almost sure if (a$'$) and (b$'$) hold. 
\endproclaim

\proclaim {Theorem 2} Assume that $0<\a<\a_c'$, and that 
(a)--(c) hold. Then, for any fixed $s$
$$
P\(H\le cm-(1-\a/\a_c')\,m\,G^{-1}(s/n)\,\mid\, p_1,\dots,p_n\)\to e^{-s}
$$
in probability. In particular, 
$$
P\(H\le cm-(1-\a/\a_c')\,m\,G^{-1}(s/n)\)\to e^{-s}. 
$$
\endproclaim

Throughout, we follow the usual convention in defining 
$G^{-1}(x)=\sup\{y: G(y)<x\}$ to be the left continuous 
inverse of $G$, although any other inverse works as well. 

Assume, for simplicity, that, as $x\to 0$, 
$G(x)$ behaves as $x^\eta$ for some $\eta>2$. Then, in 
contrast with the pure regime, the annealed fluctuations 
in composite regime  
scale as $m^{1-1/\eta}$, while the 
quenched ones scale as $m^{1/2}$. 
In fact, this can be guessed from [GTW2]. Namely, 
as explained in Section 2 of that paper, 
the maximal increasing path has a nearly vertical 
segment of length asymptotic to $(1-\a/\a_c')m$ in (or near) the 
column of $A$ which uses the largest
probability $p_1$. Therefore, this vertical part 
of the path dominates the fluctuations, as
the rest presumably has $o(\sqrt m)$ fluctuations. 
(These are most likely {\it not\/} of the order 
exactly $m^{1/3}$ as they correspond to the 
critical case $\a=\a_c'$. The precise nature 
of the critical fluctuations is an interesting 
open problem.) The variables in the $p_1$--column 
are Bernoulli with variances about $b(1-b)$, 
thus the contribution of the vertical part 
to the variance is about $(b(1-b)(1-\a/\a_c')m)^{1/2}=\tau\sqrt n$. 
The annealed case then simply picks up the 
variation in the extremal statistic $p_1$. 

Simple as the above intuition may be, Theorems 1 and 2 are 
not so easy to prove and require considerable 
additional technical details. We also note the mysterious correction 
$2\tau\sqrt n$ in Theorem 1 for which 
we have no intuitive explanation.

The fluctuations 
results in [GTW2] and the present paper
thus sharply distinguish between two 
different phases of {\it one\/}  
particular growth model. 
Nevertheless, it seems natural to speculate  
that this phenomenon is universal in the sense 
that it occurs in other one--dimensional 
finite range dynamics of ODB type, started from 
a variety of initial states.
Indeed, such universality has been established in other
random matrix contexts [Sos]. 
Fluctuations of higher--dimensional versions 
seem much more elusive; 
it appears that a glassy transition should take place, 
but the fluctuation 
scalings could be completely different. 
 
To elucidate, we present 
some simulation results. 
In all of them, we 
start from the flat substrate $h_0\equiv0$ and 
use $F(s)=1-(1-2s)^\eta$, so that $b=1/2$. It is expected that, as
$\eta$ increases, the quenched fluctuation 
experience a sudden jump from $1/3$ to $1/2$. 
We simulate two dynamics, the ODB and the two--sided
digital boiling (abbreviated simply as DB), given 
by 
$$
h_{t+1}(x)=\max\{h_{t}(x-1), h_t(x+1), h_{t}(x)+\e_{x,t}\}. 
$$
The top of Figure 1 illustrates the ODB on 600 sites (with 
periodic boundary), run until 
time 600. The occupied sites are periodically colored so that the sites 
which become occupied at the same time are given the 
same color. On the left, $\eta=1$ (i.e., 
$p$ is uniform on $[0,1/2]$ and $\a_c'=0$),
while $\eta=3$ (and $\a_c'>0$) on the right. 
The darkly colored sites thus give the height of the 
surface at different times and provide a glimpse of 
its evolution. In the pure regime ($\eta=1$), the boundary 
of the growing set reaches a local equilibrium 
([SK], [BFL]), while in the composite regime ($\eta=3$)
the boundary apparently divides into domains, which are
populated by different equilibria and grow sublinearly.
This is the mechanism that causes increasing fluctuations. 
The bottom of Figure 1 confirms this observation; it features 
a log--log plot of quenched standard deviation (estimated
over 1000 independent trials) of $h_t(0)$ vs. $t$ up to $t=10\,000$.
The $\eta=1$ case is drawn with $+$'s  
and the $\eta=3$ 
case with $\times$'s; 
the two least squares approximations lines 
(with slopes 0.339 and 0.517, respectively) are also drawn. 
We note that the asymptotic speed of this
flat interface is known: $\lim_{t\to\infty} h_t(0)/t=\sup_{\a>0}
(\a+1)c(\a)$.  Here is the reason: if ODB dynamics $h_t^i$, $h_t$  start 
from initial states $h_0^i$, $h_0=\sup_i h_0^i$, respectively, 
and are coupled by using the same coin flips $\e_{x,t}$, 
then $h_t=\sup_i h_t^i$ for every $t$.

\vskip1cm
\epsfysize=0.1235\vsize
\centerline{\epsffile{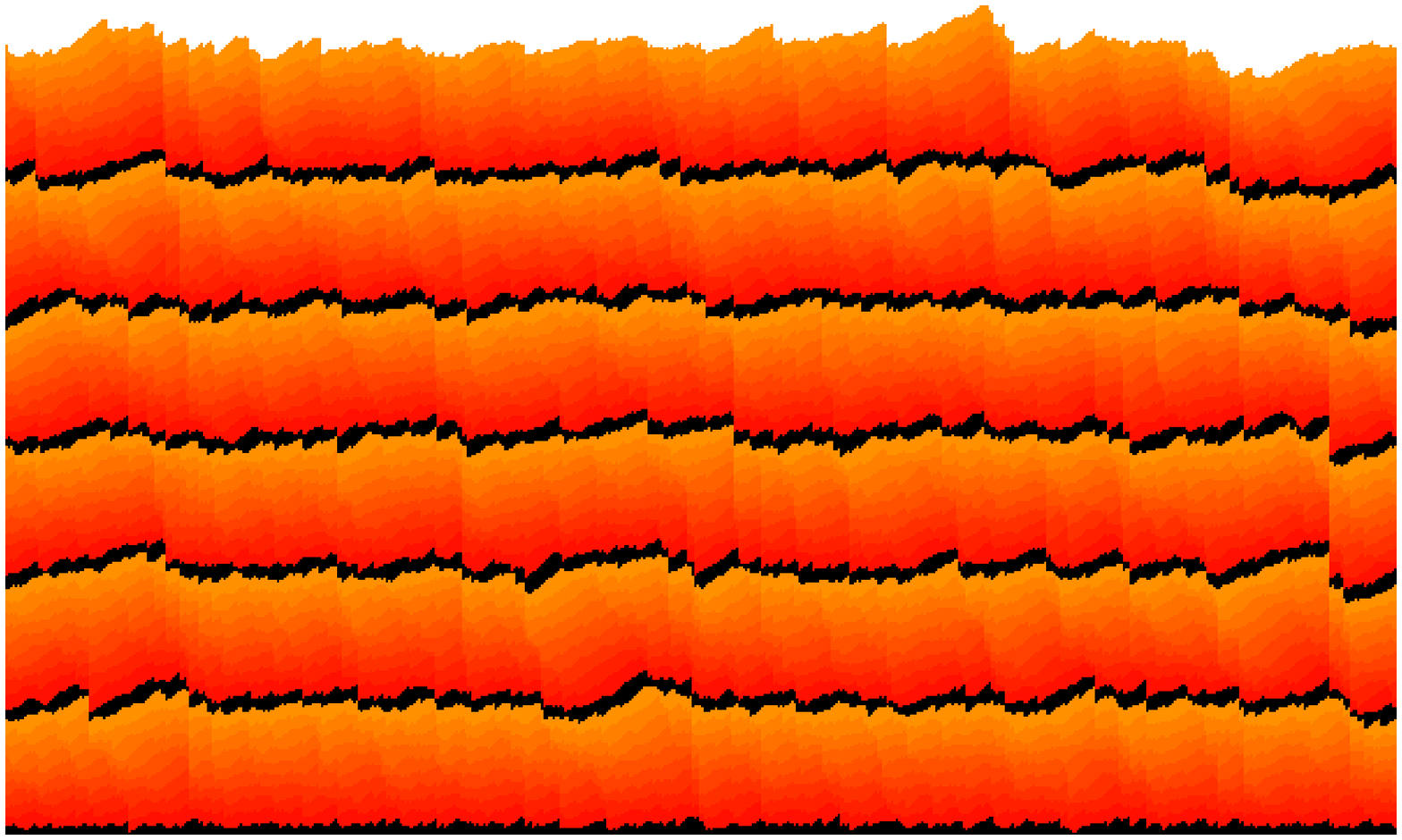}\hskip0.2cm\epsfysize=0.1235\vsize
\epsffile{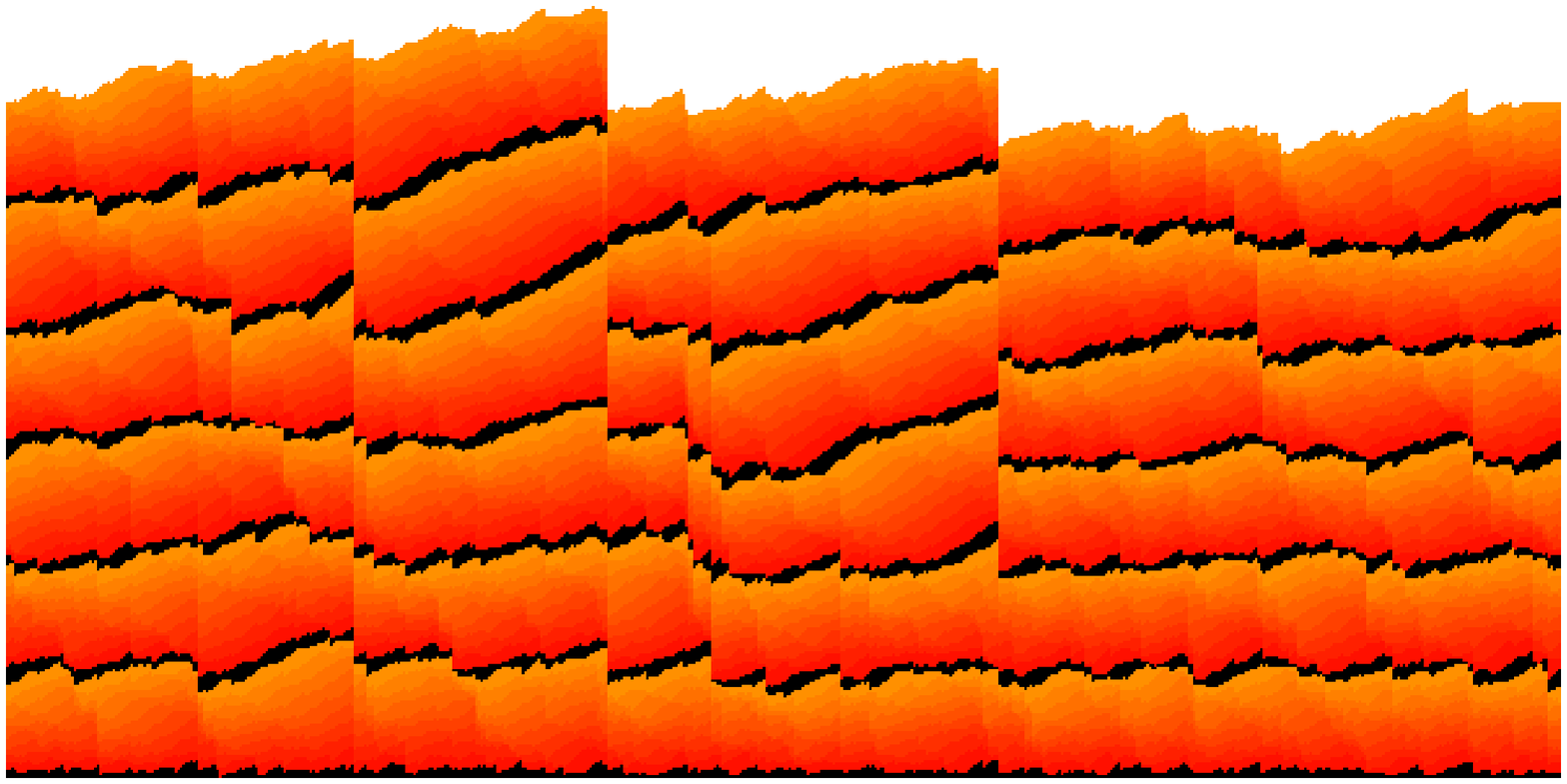}}
\smallskip
\epsfysize=0.30\vsize
\centerline{\epsffile{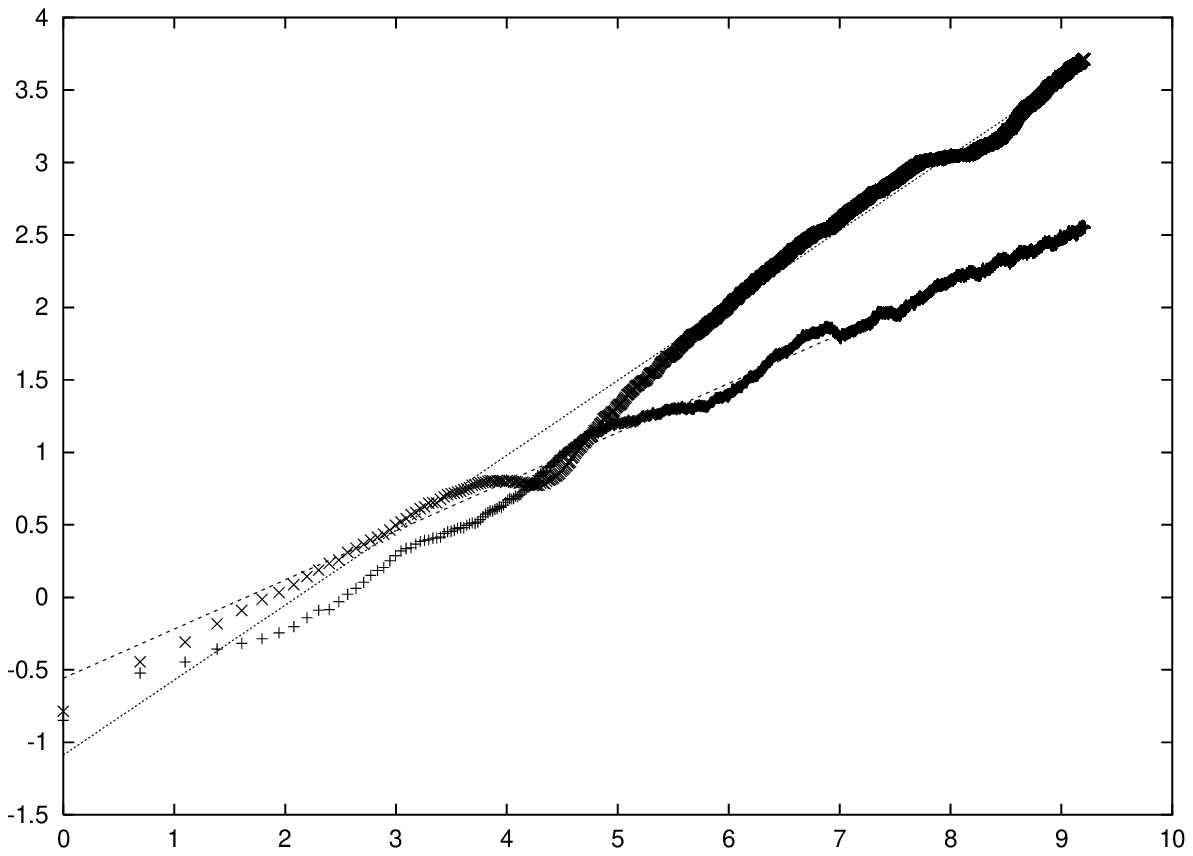}} 
\smallskip
\centerline {Figure 1. Evolution and quenched deviation 
in the two phases of disordered ODB.}
\smallskip

Perhaps surprisingly, it appears that the phase transition 
in the DB does {\it not\/} occur at $\eta=2$, and 
in general the delineation is much murkier. At this point, 
we cannot even eliminate the possibility of continuous 
dependence of fluctuation exponent on $\eta$. 
In Figure 2, we present 
the results of simulations for $\eta=0.2$ (left) and $\eta=1$
(right). The top figures only show 
evolution near time $t=5000$, as no difference is 
readily apparent at earlier times. 
The plot of quenched deviations is analogous to 
the one in Figure 1, with the least squares 
slopes 0.395 ($\eta=0.2$) and 0.49 ($\eta=1$).

\epsfysize=0.19\vsize
\centerline{\epsffile{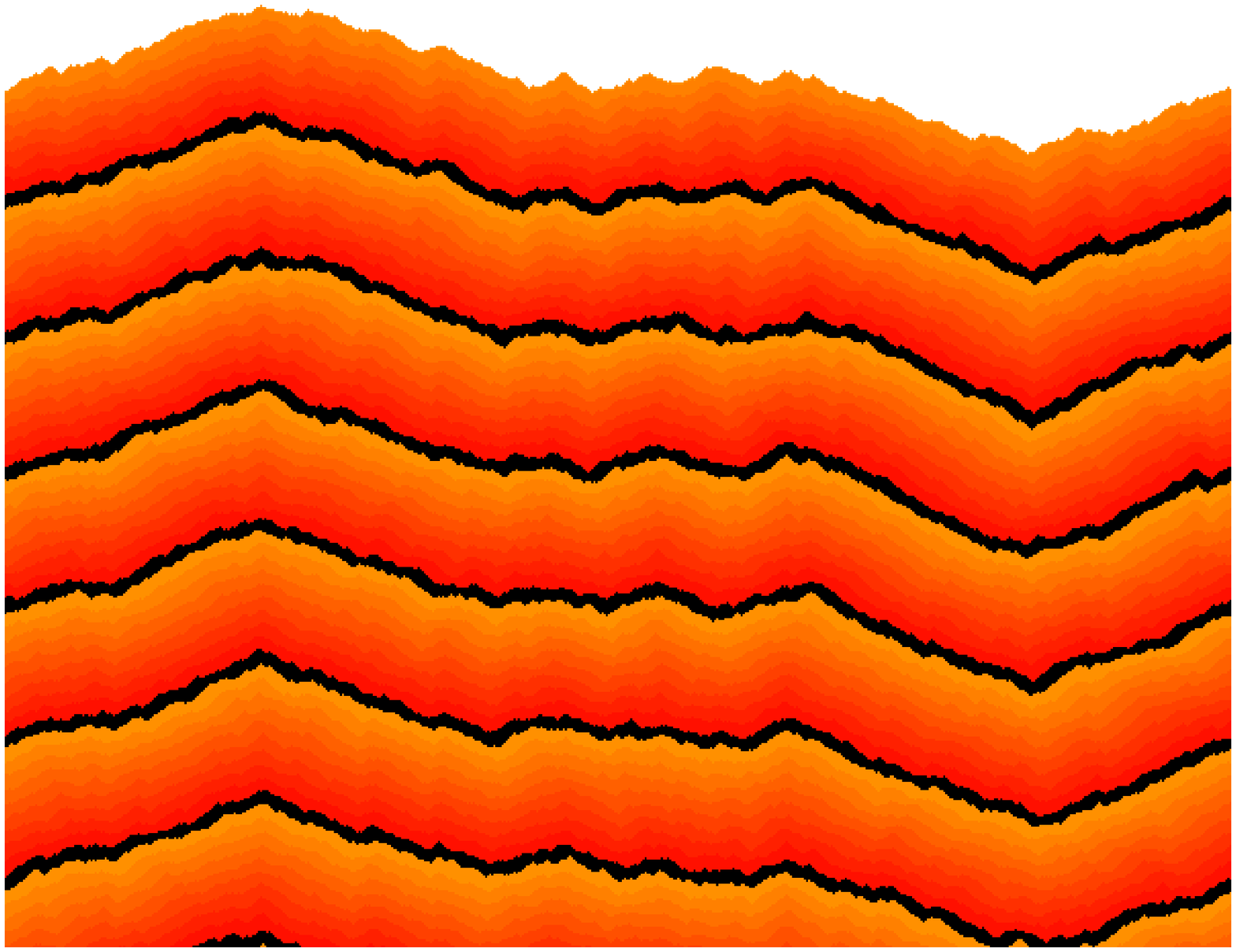}\hskip0.2cm\epsfysize=0.19\vsize
\epsffile{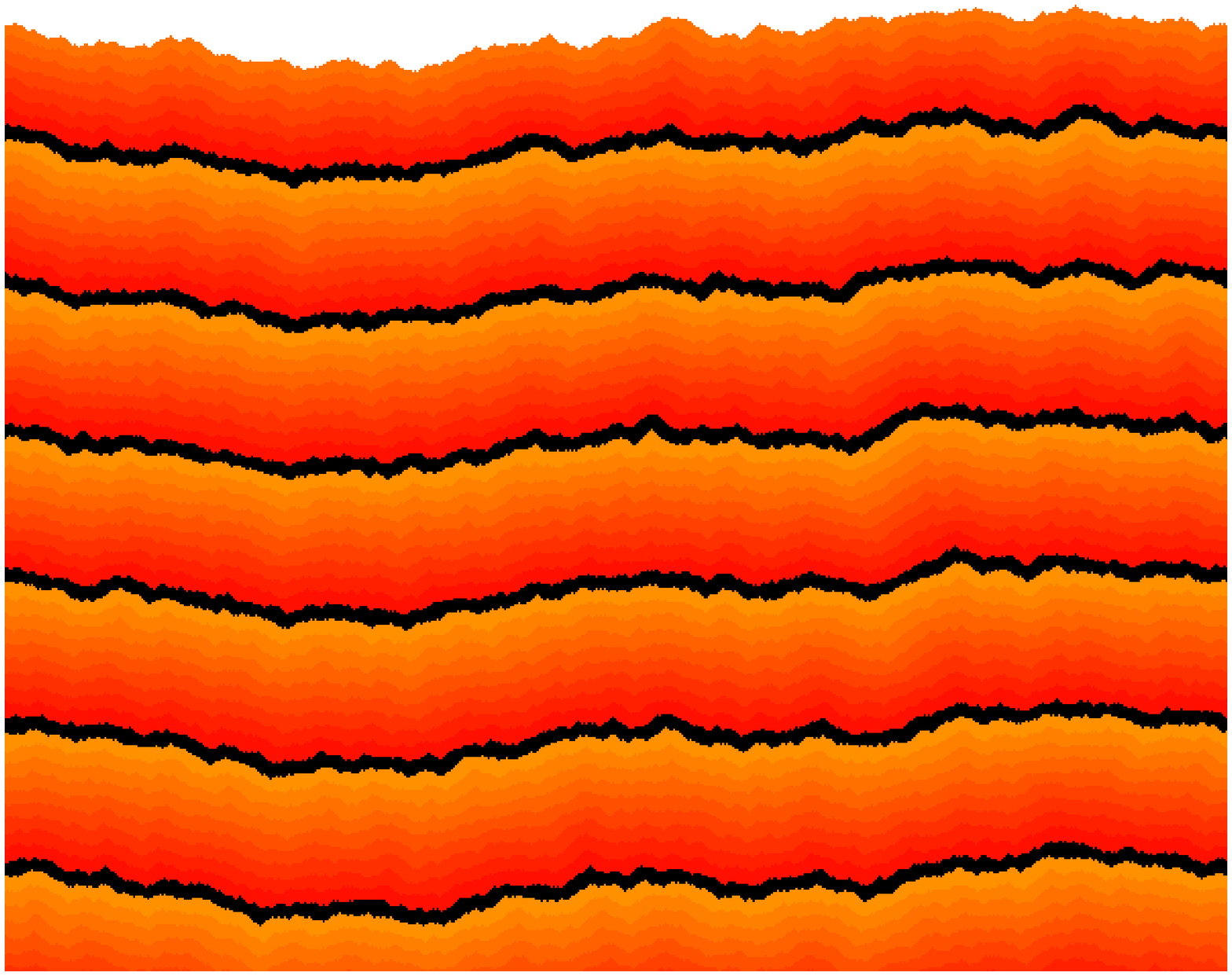}}
\smallskip
\epsfysize=0.30\vsize
\centerline{\epsffile{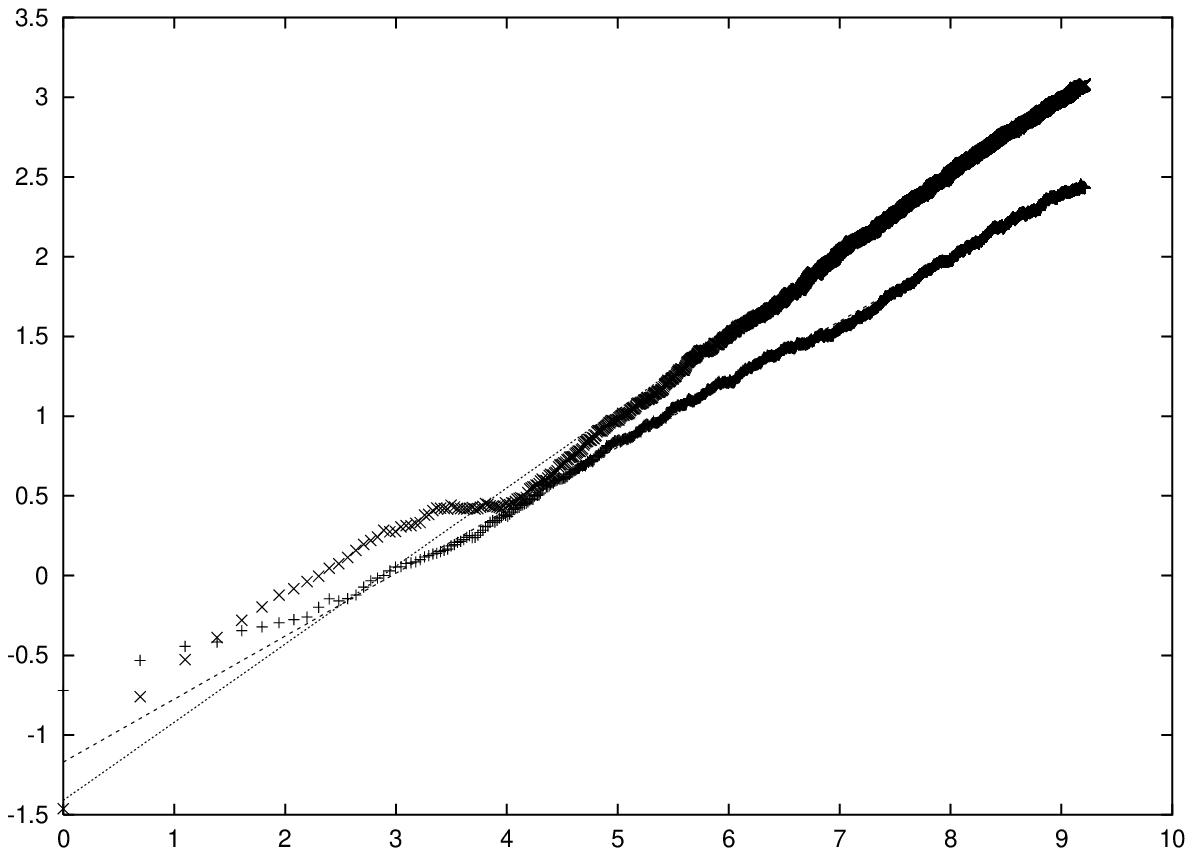}} 
\smallskip
\centerline {Figure 2. Evolution and quenched deviation 
in disordered DB.}
\smallskip

The organization of the rest of the paper is as follows. Section 2 
reviews the set-up from [GTW1, GTW2], in Section 3 we prove
the relevant asymptotic properties of the order statistics and  
of the solutions of (1.1) and (1.2), and demonstrate how 
Theorem 2 follows from Theorem 1.  
Section 4 is a detailed analysis 
of the asymptotic behavior of steepest descent curves. 
The proof of convergence in probability in Theorem 1 is then 
concluded in Section 4. Finally, Section 5 strengthens the 
results of Section 3 (under the stronger conditions) so that 
almost sure convergence is implied. 
 
\vfill\eject
\subheading {2. The basic set-up}

\define \cb{\bar c} 
\define \Ph{\phi} 
\define \ga{\gamma}

We recall how we approached these problems in [GTW1,GTW2]. The starting point
is the identity
$$\Pr(H\le h)=\det\,(I-K_h),$$
where $K_h$ is the infinite matrix acting on $\ell^2({\bZ}^+)$ with
$(j,k)$--entry
$$K_h(j,k)=\sum_{\l=0}^{\iy}(\ph_-/\ph_+)_{h+j+\l+1}\;
(\ph_+/\ph_-)_{-h-k-\l-1}.$$
The subscripts denote Fourier coefficients and the functions
$\ph_{\pm}$ are given by
$$\ph_+(z)=\prod_{j=1}^n(1+r_jz),\ \ \ \ph_-(z)=(1-z\inv)^{-m}.$$
The matrix $K_h$ is the product of two matrices, with $(j,k)$--entries given
by
$$
\align 
&(\ph_+/\ph_-)_{-h-j-k-1}={1\ov 2\pi i}
\int\prod_{j=1}^n(1+r_jz)\;(z-1)^m\,z^{-m+h+j+k}\,dz,\\
&(\ph_-/\ph_+)_{h+j+k+1}={1\ov 2\pi i}
\int\prod_{j=1}^n(1+r_jz)\inv\;(z-1)^{-m}\,z^{m-h-j-k-2}\,dz.
\endalign
$$
The contours for both integrals go around the origin once
counterclockwise; in the second
integral 1 is on the inside and all the $-r_j\inv$ are on the outside.

If $h=c_n\,m+h'$ we have
$$
\align 
&(\ph_+/\ph_-)_{-h-j-k-1}={1\ov 2\pi i}
\int\ps(z)\,z^{h'+j+k}\,dz,\tag{2.1}\\
&(\ph_-/\ph_+)_{h+j+k+1}={1\ov 2\pi i}\int\ps(z)\inv\,z^{-h'-j-k-2}\,dz,
\tag{2.2}
\endalign
$$
where
$$\ps(z)=\prod_{j=1}^n(1+r_jz)\;(z-1)^m\,z^{-(1-c_n)\,m}.$$
The idea is to apply steepest descent to the above
integrals. If $\si(z)=m\inv\log\,\ps(z)$ then
$$\si'(z)={\al\ov n}\,\sum_{j=1}^n{r_j\ov 1+r_jz}+{1\ov z-1}+{c_n-1\ov 
z}\tag{2.3}$$
and, with $u_n$ and $c_n$ as defined above, $\si'(u_n)=\si''(u_n)=0$. The 
steepest descent curves both
pass through $u_n$. As $n\ra\iy$ the zeros/poles $-r_j\inv$ accumulate
on the half-line $(-\iy,\,\xi]$ where $\xi=1-b\inv$. In the pure regime the 
points $u_n$
and the curves are bounded away from this half-line, behave regularly and have 
nice
limits. However in the composite regime the points and curves come very close to
$\xi$, their behavior is not so simple, and we apply steepest descent not quite 
as
described. 

\subheading{3. Preliminary lemmas I: properties of $p_n$, $u_n$, and $c_n$}

Until Section 5, we assume that all limits are in probability, 
unless otherwise indicated. 
To prove the first part 
of Theorem 1 and Theorem 2, we thus assume that 
(a)--(c) hold.  

We let $q_j=b-p_j$, so that $q_1,\,\cd,q_n$ are chosen 
independently according to the
distribution function $G$, then ordered so that $q_1\le q_2\le\cd\le q_n$.

Let $t_1< t_2< \dots < t_n$
be an ordered sample of i.i.d.\ uniform $(0,1)$ 
random variables.
Then we may construct the $G$--sample 
by setting $q_j=G^{-1}(t_j)$. We will also use the well-known fact
that,
given $t_j$, the conditional distribution
of $t_1,\dots t_{j-1}$ is that of an ordered sample of 
$j-1$ uniforms on $[0,t_j]$.

\proclaim{Lemma 3.1} There exist a positive constant 
$c_1$ 
so that $x\le G(G^{-1}(x))\le x/c_1$ for $x\in (0,1)$. 
Moreover,
$G(G^{-1}(x))\sim x$
as $x\to 0$. 
\endproclaim

\demo{Proof} 
Write the complement of the range of $G$ as $\cup_i I_i$, 
where $I_i$ are disjoint and either 
of the form $[a_i,b_i)$ or $(a_i, b_i)$.
If $x\in (0,1)$ is in the range of $G$, then 
$G(G^{-1}(x))=x$, otherwise, if $x\in I_i$,
$G(G^{-1}(x))=b_i$. By (a), $b_i\sim a_i$ if
$a_i\to 0$. The last sentence in the statement 
is then proved, and the first follows. $\square$
\enddemo
 
\proclaim{Lemma 3.2} With $c_1$ as in Lemma 3.1, 
for $\eta<1$ and $j\ge 2$,
$$
\Pr\left(G(q_1)>\eta G(q_j)\right)\le
(1-c_1\eta)^{j-1}.
$$
\endproclaim

\demo{Proof} By Lemma 3.1 and remarks preceding it,  
$$
\Pr\(G(q_1)>\eta G(q_j)\)
\le \Pr\(t_1>{c_1}\eta t_j\)=
\(1-{c_1}\eta\)^{j-1}.\qquad\square
$$
\enddemo

\proclaim{Lemma 3.3} $\lim_{n\ra\iy}\Pr\left(q_1\le 
G\inv(s/n)\right)=1-e^{-s}$.
\endproclaim

\demo{Proof} Fix an $\e>0$.  First, by
monotonicity of $G^{-1}$, 
$t_1\le s/n$ implies $q_1\le G^{-1}(s/n)$.
Second, by Lemma 3.1 and the monotonicity of $G$ we have that, 
for large enough $n$,
$q_1\le G^{-1}(s/n)$ implies $t_1\le
G(G^{-1}(t_1))=G(q_1)\le G(G^{-1}(s/n)
\le (1+\e)s/n$. These give the inequalities
$P(q_1\le G^{-1}(s/n))\ge 1-(1-s/n)^n$, and 
$P(q_1\le G^{-1}(s/n))\le 1-(1-(1+\e)s/n)^n$. 
The statement of the lemma now follows upon first letting $n\to\infty$
and then 
$\e\to 0$.
 $\square$
\enddemo

\comment 

\demo{Proof} Let $t_1$ be as in the previous proof.
Fix an $\e>0$. Find a $\d>0$ so that $x<\d$ 
implies $x\le G(G\inv(x))<(1+\e)x$.
Then $P(G(G\inv(t_1))<s/n)\le
P(t_1< s/n)+P(t_1>\d)\to 1-\exp(s)$.
Together with a similarly proved lower bound this implies
that $P(q_1<G\inv(s/n))=P(G(q_1)<s/n)\to 1-\exp(-s)$.
Thus we only need to show that 
$P(G^{-1}(t_1)=G^{-1}(s/n))=
P(t_1$ and $s/n$ are in the same $I_i)\to 0$.
(Here $I_i$ are in the same as in the proof 
of Lemma 3.1.) But the probability in question 
is either 0 if $s/n$ is in the range of $G$ or,
if $s/n\in I_i$,  is bounded 
above by $P(t_1>a_i\,|\,t_1<b_i)\le (b_i-a_i)/b_i$. 
However, $a_i\to 0$ as $n\to\infty$,
thus $a_i\sim b_i$ 
and the proof is concluded. $\square$
\enddemo 

\endcomment

\noi{\bf Remark}. It follows from Lemma 3.3,
and the fact that $G(x)=o(x^2)$ near $x=0$,
that $n^{1/2}q_1\ra\iy$ as $n\ra\iy$.\sp

\proclaim{Lemma 3.4} 
With high probability $q_1/q_2$ is bounded away from 1 as 
$n\ra\iy$.
More precisely, for every $\eta>0$ there is a $\dl>0$ such that
$\Pr(q_1\le(1-\dl)\,q_2)\ge 1-\eta$ for large enough $n$.
\endproclaim

\demo{Proof} 
It follows from Lemma 3.1 
that for every $\eta>0$ there exists a $\d_1>0$ so that
the following implication holds for $t_2<\d_1$:
if $G(q_1)>(1-\d_1)G(q_2))$ then $t_1>(1-\eta)t_2$.
Furthermore, by the assumption (a), there exists a
$\d\in (0,\d_1)$ so that, for $t_2<\d$, $q_1>(1-\d)q_2$ implies
$G(q_1)>(1-\d_1)G(q_2)$. Therefore, 
$$
P(q_1>(1-\d)q_2)\le P(t_1>(1-\eta)t_2)+P(t_2>\d)=\eta +P(t_2>\d), 
$$
and the proof is concluded since $t_2\to 0$ a.s.
$\square$
\enddemo

\proclaim{Lemma 3.5}
$n\inv\sum_1^nq_1/q_j^3\ra 0$ as $n\ra\iy$. 
\endproclaim

\demo{Proof}
For any fixed $k$ we have $n\inv\sum_{j=1}^kq_1/q_j^3\le 
k/nq_1^2\to 0$.
Also, $n\inv\sum_{j=k+1}^n q_j^{-2}<\lan q^{-2}\ran+1$ a.s. 
for large $n$.

Let $\dl>0$ be given.  By the above paragraph, 
it suffices to show that
$$\limsup_{n\to\infty}\Pr\left({q_1\ov q_{k+1}}> \dl \right)$$
will be arbitrarily small for sufficiently large $k$.
Now, from the assumption (b), 
it follows that for some $\eta>0$ we have 
$G(q_1)> \eta G(q_{k+1})$ 
whenever
$q_1> \dl q_{k+1}$ and $q_1<\eta$. 
With this $\eta$ (which we may assume is less than 1) 
we have, from Lemma 3.2,
$$\Pr\left({q_1\ov q_{k+1}}>\dl\right)
\le(1-c_1\eta)^{k}+P(q_1\ge \eta),$$
which is clearly enough.  
$\square$
\enddemo

From now on $\{\ph_n\}$ will denote a sequence 
of random variables satisfying $\ph_n=o(q_1)$. 
Since $q_1\gg n^{-1/2}$ we shall assume when convenient that also $\ph_n\gg 
n^{-1/2}$. In the statement of the next lemma, the expression $O(\ph_n)$ 
could have been replaced by the less awkward $o(q_1)$.
The reasons for the present statement are that
the substitute for this lemma (Lemma 6.2) when we consider 
almost sure convergence will have this form, and that  
the same sequence $\{\ph_n\}$ will
appear in later lemmas.

\proclaim{Lemma 3.6}
Let $\{v_n\}$ be a sequence of points in a disc with 
diameter
the real interval $[-r_1\inv-O(\ph_n),\,\xi]$. Then
$$\lim_{n\ra\iy}{1\ov 
n}\sum_{j=2}^n{r_j\ov(1+r_jv_n)^2}=\lan{r\ov(1+r\xi)^2}\ran.$$
\endproclaim

\demo{Proof} Write $v_n=(b_n-1)/b_n$. Then if we recall that $\xi=(b-1)/b$ 
and
$p_j=b-q_j$ we see that $b-b_n$ lies in a disc with diameter $[0,\,q_1+O(\ph_n)]$ 
and that
$${1\ov 
n}\sum_{j=2}^n{r_j\ov(1+r_jv_n)^2}={1\ov n}\sum_{j=2}^n{b_n^2(b-q_j)(1-b+q_j)\ov 
(b_n-b+q_j)^2}.$$
If we subtract from this the same expression with $b_n$ replaced by $b$,
that is, 
$${1\ov n}\sum_{j=2}^n{b^2(b-q_j)(1-b+q_j)\ov q_j^2}, \tag{3.1}$$
 we obtain
$${1\ov n}\sum_{j=2}^n(b-q_j)(1-b+q_j)
\left[{b_n^2\ov(b_n-b+q_j)^2}-{b^2\ov q_j^2}\right].\tag{3.2}$$
We shall show that this is $o(1)$. Assuming this for the moment,
we can finish the proof by first noting that we may, with error o(1), 
start the sum in (3.1) at $n=1$ since $q_i\gg n^{-1/2}$, and then 
(3.1) has the 
a.s. limit
$$\lan {b^2(b-q)(1-b+q)\ov q^2}\ran=\lan{r\ov(1+r\xi)^2}\ran.$$

It remains to show that (3.2) is $o(1)$.
If we replace the numerator $b^2$ on the right by $b_n^2$, the error is
$o(1)$, since $n\inv\sum q_j^{-2}$ is a.s. bounded. If we make this replacement 
then what
we obtain is bounded by a constant times
$${b\ov n}\sum_{j=2}^n\left|{(b_n-b)^2-2(b_n-b)q_j\ov 
q_j^2(b_n-b+q_j)^2}\right|.$$
Since  
$|b-b_n|\le q_1+O(\ph_n)=q_1+o(q_1)$ it follows from Lemma~3.4 that 
$|b_n-b+q_j|$ 
is at least a constant times $q_j$ for large $n$ and so the
above is at most a constant times
$${1\ov n}\sum_{j=2}^n{|b_n-b|\ov q_j^3}\le
{1\ov n}\sum_{j=2}^n{q_1\ov q_j^3},$$
and by Lemma 3.5 this is $o(1)$. $\square$
\enddemo

We denote 
$$
\theta=1-\a/\a_c', \qquad
\b=\left({(1-b)\,\al\ov b^3\,\theta}\right)^{1/2}.
\tag 3.3
$$

\proclaim{Lemma 3.7} We have $u_n=-r_1\inv+\b n^{-1/2}+o(n^{-1/2})$ 
as $n\ra\iy$. 
\endproclaim
 
\demo{Proof} We show first that $u_n\ge\xi$ cannot occur for arbitrarily large $n$.
If it did, then we would have, using equation (1.1) for $u_n$, 
$$b^2={1\ov (\xi-1)^2}\le{1\ov (u_n-1)^2}
\le {\al\ov n}\,{r_1\ov(1+r_1\xi)^2}+{\al\ov n}\,\sum_{j=2}^n{r_j\ov
(1+r_j\xi)^2}.$$
It follows from the remark following Lemma 3.3 that the first term on the right
is $o(1)$ and from Lemma~3.6 that the second term on the right has limit
$$\al\lan {r\ov(1+r\xi)^2}\ran=\al b^2\lan{p(1-p)\ov(b-p)^2}\ran<b^2$$
since we are in the composite regime. This contradiction shows that $u_n\le\xi$
for sufficiently large $n$, and so
$u_n\in[-r_1\inv,\,\xi]$. By Lemma~3.6 again, 
$${\al\ov n}\,\sum_{j=2}^n{r_j\ov
(1+r_ju)^2}={1\ov(u-1)^2}
\to\al\lan {r\ov(1+r\xi)^2}\ran=b^2\al/\al_c'.$$
Therefore the equation (1.1) for $u_n$ becomes
$${\al\ov n}\,{r_1\ov(1+r_1u_n)^2}={1\ov (\xi-1)^2}-\al\lan 
{r\ov(1+r\xi)^2}\ran+o(1)
=b^2\theta+o(1).$$
Since $r_1=b/(1-b)+o(1)$ we find that the solution is as stated. $\square$
\enddemo

Next, we see how $c_n$ behaves. 

\proclaim{Lemma 3.8} We have $c_n=c(\al,F)-\theta \,q_1+o(q_1)$ as $n\ra\iy$,
where $\theta$ is given in (3.3).
\endproclaim

\demo{Proof} Write
$$c_n={1\ov 1-u_n}-{\al\ov n}\sum_{j=2}^n{r_ju_n\ov 1+r_ju_n}-
{\al\ov n}{r_1u_n\ov 1+r_1u_n}.\tag{3.4}$$
By Lemma 3.7, the last term above is $O(n^{-1/2})$. Equation (1.1)
tells us that
$${d\ov du}\left({1\ov 1-u}-{\al\ov n}\sum_{j=1}^n{r_ju\ov 
1+r_ju}\right)\Big|_{u=u_n}=0,$$
and so
$${d\ov du}\left({1\ov 1-u}-{\al\ov n}\sum_{j=2}^n{r_ju\ov 
1+r_ju}\right)\Big|_{u=u_n}=
{\al\ov n}{r_1\ov (1+r_1u_n)^2}={\al\ov r_1\b^2}+o(1)={\al(1-b)\ov b\b^2}+o(1).$$

By Lemma 3.6 and its proof, 
with an error $o(1)$ the derivative of the expression in the 
parentheses above equals in
$[u_n,\xi]$ what it equals at $u=\xi$, so the above holds with $u_n$ replaced 
by any
point in this interval. From this and (3.4) we get
$$c_n=c(u_n)=c(\xi)-{\al(1-b)\ov b\b^2}(\xi-u_n)+o(\xi-u_n).$$
We have
$$\xi-u_n=1-b\inv-r_1\inv+O(n^{-1/2})=p_1\inv-b\inv+O(n^{-1/2})
={q_1\ov b^2}+o(q_1),$$
where we have used the fact that $q_1\gg n^{-1/2}$. Thus
$$ c_n=c(\xi)-{\al(1-b)\ov b^3\b^2}q_1+o(q_1).$$
Finally, as $\<(b-p)^2\><\infty$, we can 
use the central limit theorem  
to conclude that $c(\xi)=c(\al,F)+O(n^{-1/2})$,
which completes the proof. $\square$
\enddemo

\noi{\bf Remark}. Lemmas 3.3 and 3.8 show 
that Theorem 2 follows from the part of Theorem 1
on convergence in probability.\sp

\subheading{4. Preliminary lemmas II: steepest descent curves} 

\def\AH {4.1}

\def\siders {4.2}
\def\sonezz {4.3}

\def\BH {4.4}

\def\sioneest {4.5}
\def\allcn {4.6}
\def\siexp {4.7}
\def\uplus {4.8}
\def\siuplus {4.9}
\def\uest {4.10}
\def\logu {4.11}
\def\logeq {4.12}

\def\CH {4.13}

\def\DH {4.14}

\def\regionone {4.15}
\def\regiontwo {4.16}
\def\kernel {5.1}

Now we go to our integrals (2.1) and (2.2). We are not going to apply steepest
descent with $\ps$ as the main integrand, but rather with the function $\ps_1$
which is $\ps$ with
the factor $1+r_1z$ removed. It is convenient to introduce the notation
$$\ps_1(z,c)=\prod_{j=2}^n(1+r_jz)\;(z-1)^m\,z^{-(1-c)\,m},$$
where $c>0$. (This parameter is not to be confused
with the time constant $c=c(\al,F)$ defined earlier.) Thus $\ps_1(z)=\ps_1(z,c_n)$ in
this notation. We also define the integrals
$$I^+(c)={1\ov 2\pi i}\int(1+r_1z)\,\ps_1(z,c)\,dz,
\ \ \ I^-(c)={1\ov 2\pi i}\int(1+r_1z)\inv\,\ps_1(z,c)\inv z^{-2}\,dz.$$
(Since $I^+(c)=0$ when $c\ge1$ we always assume that $c<1$.)
Notice that these are exactly the integrals (2.1) and (2.2) when we set
$$c=c_n+(h'+j+k)/m.$$
Since $j,k\ge0$ and we will eventually set $h'=sn^{1/2}$,
we may also assume that 
$$c\ge c_n-O(n^{-1/2}).\tag{\AH}$$

To apply steepest descent to $I^{\pm}(c)$ we must locate the critical points and
determine the critical values of $\ps_1(z,c)$. Thus we define
$$\si_1(z,c)={1\ov m}\log\,\ps_1(z,c),$$
so that 
$$\si_1'(z,c)={\al\ov n}\,\sum_{j=2}^n{r_j\ov 1+r_jz}+{1\ov z-1}+{c-1\ov z}.$$
As before, if the parameter $c$ does not appear we take it to be $c_n$, 
e.g., $\si_1(z)=\si_1(z,c_n)$. So
$$\si_1'(z)={1\ov m}\log\,\ps_1(z)=\si'(z)-{\al\ov n}{r_1\ov 1+r_1z}.$$
Using $\si'(u_n)=\si''(u_n)=0$ we get from
the above and Lemma~3.7 that
$$\si_1'(u_n)=-{\al\ov\b\sqrt n}(1+o(1)),\ \ \ \si_1''(u_n)={\al\ov\b^2}(1+o(1)).\tag{\siders}$$

To determine the critical values of $\si_1(z,c)$ let us first find
the value of $c$ for which its derivative has a double zero. (This is the
analogue of the quantity $c_n$ for $\si(z)$.) For this we use the analogue of
(1.1) and (1.2) but where the terms corresponding to $j=1$ are dropped
from the sums. If we call the solution of (1.1) $\ub$ and set $\cb=c(\ub)$ then
$\si_1'(z,\cb)$ has a double zero at $\ub$.
In analogy with $u_n$, we
know that $\ub$ is to the right of and within $O(n^{-1/2})$ of $-r_2\inv$.
As for $\cb$, we use Lemma~3.8, its analogue where the sums in (1.1) 
and (1.2) start with $j=2$,
as well as Lemma~3.4, to see that to a first approximation
$$\cb=c_n-\theta(q_2-q_1)$$
and that $q_2-q_1\gg n^{-1/2}$. From this and (\AH) we see that $c>\cb$.

Using subscripts for derivatives now, we have
$$\si_{1z}(\ub,\cb)=\si_{1zz}(\ub,\cb)=0$$
and we want to see how the critical points $u_{c}^{\pm}$ of $\si_1(z,c)$ move away
from $\ub$ as $c$ increases from $\cb$. (Here we take
$u_{c}^-<u_{c}^+$.) The function $\si_{1z}(z,\cb)$ vanishes at $\ub$
and is otherwise positive in $(-r_2\inv,0)$. It follows that for $c$ close to but larger 
than $\cb$ we have $u_c^-<\ub<u_c^+$. Differentiating
$\si_{1z}(u_{c}^{\pm},c)=0$
with respect to $c$ gives
$$0=\si_{1zz}(u_{c}^{\pm},c)\,{du_{c}^{\pm}\ov dc}+\si_{1zc}(u_{c}^{\pm},c)
=\si_{1zz}(u_{c}^{\pm},c)\,{du_{c}^{\pm}\ov dc}+{1\ov 
u_{c}^{\pm}}.\tag{\sonezz}$$
Since $u_{c}^{\pm}<0$ it follows that
$du_{c}^+/dc\ne0$, and so each of $u_c^{\pm}$ is either a decreasing or
increasing function of $c$ for $c>\cb$. {}From their behavior that we already know for
$c$ close to $\cb$ we deduce that $u_{c}^+$ increases and $u_{c}^-$ decreases as
$c$ increases. In particular, $u_{c}^-$ is even closer to $-r_2\inv$ than $\ub$.

We remark that from (\sonezz) and the signs of $du_{c}^+/dc$ we deduce
$$\si_{1zz}(u_c^+,c)>0,\ \ \ \si_{1zz}(u_c^-,c)<0.\tag{\BH}$$

Next we shall determine the asymptotics of the critical values 
$\si(u_c^{\pm},c)$. 
The sequence $\{\ph_n\}$ is as described before Lemma 3.6.

\proclaim{Lemma 4.1} For $c-c_n=O(\ph_n)$
$$\si_1(u_{c}^+,c)=\si_1(-r_1\inv,c)-
{r_1 \b^2\ov2\al}\left(c-c_n+{2\al\ov 
r_1\b}(1+o(1))n^{-1/2}\right)^2\tag{\sioneest}$$
and for all $c\ge c_n$
$$\si_1(u_{c}^+,c)< \si_1(-r_1\inv,c)-\eta 
n^{-1/2}\,(c-c_n)+O(n^{-1}).\tag{\allcn}$$
for some $\eta>0$. Moreover for all $c$
$$\si_1(u_{c}^-,c)>\si_1(-r_1\inv,c)+\ph_n^2$$
when $n$ is sufficiently large.
\endproclaim

\noi{\bf Remark}. In these and analogous inequalities below we think of $\si_1$ 
as actually
meaning $\Re\si_1$.

\demo{Proof} Consider first the case $c=c_n$. We have
$$\si_1(u_n+\z)=\si_1(u_n)+\si_1'(u_n)\,\z+\z^2\int_0^1(1-t)\,\si_1''(u_n+t\z)\,d
t.$$
If $\z=O(\ph_n)$ then it follows from Lemma~3.6 that 
$\si_1''(u_n+t\z)=\si_1''(u_n)+o(1).$
Hence, by ({\siders}), we have for such $\z$
$$\si_1(u_n+\z)=\si_1(u_n)-{\al\ov\b\sqrt 
n}\z+\left({\al\ov2\b^2}+o(1)\right)\z^2.
\tag{\siexp}$$
This has zero derivative for
$$\z={\b\ov \sqrt n}(1+o(1))$$
and it follows that
$$u_{c_n}^+=u_n+{\b\ov \sqrt n}(1+o(1))=-r_1\inv+{2\b\ov \sqrt 
n}(1+o(1)).\tag{\uplus}$$
(This critical value must be $u_{c_n}^+$ rather than $u_{c_n}^-$ since the latter 
is within
$O(n^{-1/2})$ of $-r_2\inv$.) From this and ({\siexp}), taking $\z=-r_1\inv-u_n
=-(\b+o(1))n^{-1/2}$ and $\z=u_{c_n}^+-u_n=(\b+o(1))n^{-1/2}$ and subtracting,
it follows that
$$\si_1(u_{c_n}^+)=\si_1(-r_1\inv)-2(\al+o(1))n^{-1}.\tag{\siuplus}$$

To determine the behavior of $u_{c}^+$ and $\si_1(u_{c}^+,c)$ for more general 
$c$
we assume first that
$$c=c_n+o(1),\ \ \  u_c^+=u_n+O(\ph_n)=-r_1\inv+O(\ph_n).$$
Then
$$\si_{1zz}(u_{c}^+,c)=\si_1''(u_n)-{c-c_n\ov {u_c^+}^2}={\al\ov\b^2}+o(1)$$
by ({\siders}). Therefore ({\sonezz}) gives
$${du_{c}^+\ov dc}=-(\b^2/\al+o(1))/u_{c}=r_1{\b^2\ov\al}(1+o(1),$$
whence
$$
\aligned
u_{c}^+&=u_{c_n}^++r_1{\b^2\ov\al}(c-c_n)(1+o(1))\\
&=-r_1\inv+{2\b\ov \sqrt 
n}(1+o(1))
+r_1{\b^2\ov\al}(c-c_n)(1+o(1)), 
\endaligned
\tag{\uest}
$$
by ({\uplus}). This holds if $c-c_n=O(\ph_n)$ since this assures that
$u_{c}^+=u_n+O(\ph_n)$. The above gives
$$\log (-u_{c}^+)=
\log(-r_1\inv)-2r_1\b(1+o(1))n^{-1/2}-r_1^2{\b^2\ov\al}(c-c_n)(1+o(1)).
\tag{\logu}$$
(Again, real parts are tacitly meant.)

To determine, $\si_1(u_{c}^+,c)$ we use $\si_{1z}(u_{c}^+,c)=0$
to deduce
$${d\ov dc}\si_1(u_{c}^+,c)=\log u_{c}^+.\tag{\logeq}$$
We continue to assume that $c-c_n=O(\ph_n)$ so our estimates hold.
Integrating ({\logeq}) using the first part of ({\uest}) gives
(since $u_{c_n}^+\ra-r_1\inv$)
$$
\align
\si_1(u_{c}^+,c)=&\,\,\si_1(u_{c_n}^+)+(c-c_n)\,\log u_{c_n}^+-
{1\ov2}r_1^2{\b^2\ov\al}\,(c-c_n)^2\,(1+o(1))\\
=&\,\,\si_1(-r_1\inv)-2(\al+o(1))\,n\inv+\log(-r_1\inv)(c-c_n)\\
&-2r_1\b (c-c)\,n^{-1/2}(1+o(1))
-{1\ov2}r_1^2{\b^2\ov\al}\,(c-c_n)^2\,(1+o(1)),
\endalign
$$
by ({\siuplus}) and ({\logu}). This gives ({\sioneest}).

For all $c\ge c_n$ we use the fact that
$\log(-u_{c}^+)$ is a decreasing function of $c$, since $u_{c}^+$ increases,
and integrate ({\logeq}) with respect to $c$ from $c_n$ to $c$, which gives
$$\si_1(u_{c}^+,c)\le\si_1(u_{c_n}^+)+\log(-u_{c_n}^+)(c-c_n).$$
Using ({\siuplus}) and ({\uplus}) give ({\allcn}).

For the lower bound for $\si_1(u_c^-,c)$, we assume first that $c\le c_n$.
By (\AH) this implies in particular that $c-c_n=O(n^{-1/2})$. Now $\si_1(z)$
is decreasing on the interval $(u_c^-,\,u_c^+)$ and 
$u_c^+-u_c^-\gg\ph_n$. To see the last inequality, note that,
from Lemma 3.6, $\si_{1zz}(u_n+\z,c)
\ne0$ for $\z=O(\ph_n)$ and $c-c_n=o(1)$. Therefore $\si_{1z}(u_n+\z,c)$ can 
vanish for at most
one such $\z$ and, since  $u_c^+-u_n=O(\ph_n)$, we must have 
$u_n-u_c^-\gg\ph_n$.
 
Take any sequence $\ph_n=o(q_1)$ and write
$$\si_1(u_c^-,c)\ge \si_1(u_c^+-\ph_n,c)=\si_1(u_c^+-\ph_n)+
(c-c_n)\log(\ph_n-u_{c}^+).$$
(As usual, we imagine real parts having been taken.) If we apply ({\siexp}) 
with
$\z=u_c^+-u_n$
and with $\z=u_c^+-\ph_n-u_n$ and subtract, we obtain
$$
\si(u_c^+-\ph_n)-\si(u_c^+)
={\al\ov\b} n^{-1/2}\ph_n(1+o(1))+
{\al\ov 2\b^2}\left(-2\ph_n(u_c^+-u_n)+\ph_n^2)\right)(1+o(1)).
$$
By subtracting the first parts of ({\uest}) and ({\uplus}) we see that this 
equals
$$o(n^{-1/2}\ph_n)+{\al\ov 2\b^2}\ph_n^2.$$
Since $\ph_n\gg n^{-1/2}$, as we may assume, we obtain
$$\si_1(u_c^+-\ph_n)>\si_1(u_c^+)+\eta \ph_n^2$$
for some $\eta>0$. Also, since $c-c_n>-\eta n^{-1/2}$ for some $\eta$ and
$\log(1-\ph_n/u_c^+)$ is positive and $O(\ph_n)$ we have
$$(c-c_n)\log (\ph_n-u_{c}^+)\ge (c-c_n)\log (-u_{c}^+)-\eta n^{-1/2}\ph_n.$$
Putting these together gives
$$
\si_1(u_c^-,c)>\si_1(u_c^+,c)+\eta \ph_n^2
$$
for some $\eta>0$.

This was for $c\le c_n$. For $c>c_n$ we use what we get from ({\logeq})
by replacing $^+$ with $^-$, subtracting the two, and integrating. Together
with using the already proved inequality for $c=c_n$
this gives
$$\si_1(u_c^-,c)-\si_1(u_c^+,c)>\eta \ph_n^2+\int_{c_n}^c\log(u_c^-/u_c^+)\,dc.$$
The logarithm is nonnegative. Hence
$\si_1(u_c^-,c)-\si_1(u_c^+,c)>\eta \ph_n^2$
for all $c$.

If $c-c_n=O(\ph_n)$ then using this and ({\sioneest}) give
$$\si_1(u_c^-,c)>\si_1(-r_1\inv)+\log(r_1\inv)(c-c_n)+\eta \ph_n^2.$$
with a different $\eta$. If $c\ge c_n$ we use
$$\si_1(u_c^-,c)-\si_1(u_{c_n}^-)=\int_{c_n}^c\log(-u_c^-)\,dc.$$
Since $u_c^-$ is decreasing and is less than $-r_1\inv$ when $c=c_n$ this gives
$$
\align 
\si_1(u_c^-,c)&\ge\si_1(u_{c_n}^-)+\log(r_1\inv)(c-c_n)\\
&\ge
\si_1(u_{c_n}^+)+\log(r_1\inv)(c-c_n)+\ph_n^2.
\endalign
$$
Combining this with ({\sioneest}) for $c=c_n$ shows that
$$\si_1(u_c^-,c)\ge\si_1(-r_1\inv)+\log(r_1\inv)(c-c_n)+\eta \ph_n^2$$
holds for these $c$ as well. Since $\{\ph_n\}$ was an arbitrary sequence
satisfying $\ph_n=o(q_1)$ the last statement of the lemma follows.
$\square$
\enddemo

Next we consider the steepest descent curves, which we denote by $C^{\pm}(c)$ 
corresponding to the integrals $I^{\pm}(c)$. It follows from (\BH) that $C^+(c)$
passes through $u_c^+$ because on the curve $|\ps_1(z,c)|$ has a maximum at that
point; similarly, $C^-(c)$ passes through $u_c^-$. We have enough information
to evaluate the portions of these integrals 
taken over the immediate neighborhoods
of these points, but we also have to show that the 
integrals over the rest of the curves
are negligible. This requires not only that the integrands 
are much smaller there, which they are, but also that the curves 
themselves are not too badly behaved. 

To see what is needed, let $\Ga^{\pm}$ be arcs of steepest descent 
curves for a function 
$\rho$, curves on which $\Im\rho$ is constant. 
In analogy with our $C^{\pm}(c)$
we assume $\Re\rho$ is
increasing on $\Ga^-$ as we move away from the critical point and 
decreasing on $\Ga^+$. If $s$ measures arc length on $\Ga^{\pm}$ we have for $z\in\Ga^{\pm}$
$${dz\ov ds}=\mp{|\rho'(z)|\ov\rho'(z)}.\tag{\CH}$$
If the arc goes from $a$ to $b$ then 
$$\int_{\Ga^{\pm}}|\rho'(z)|\,ds=\mp\int_{\Ga}\rho'(z)\,dz=\mp(\rho(b)-\rho(a)).$$
Hence the length of $\Ga^{\pm}$ is at most
$${|\rho(b)-\rho(a)|\ov \min_{z\in\Ga^{\pm}}|\rho'(z)|}.\tag{\DH}$$
This is to be modified if $\rho'$ has a simple
zero at $z=a$, for example. In this case we replace $\rho'(z)$ by $\rho'(z)/(z-a)$.
(This is seen by making the variable change $z=a+\sqrt{\xi}$.) 

Our goal is Lemma~4.5 below. In order to use the length estimate (\DH)
to deduce the bounds of the lemma, we must first locate regions in which our
curves are located, and then find lower bounds for $\si_1'(z,c)$
in these regions. (Upper bounds for $|\si_1(z,c)|$ will be easy.) These will be established 
in the next lemmas.

For $r>0$ define $n(r)=\#\{j:r_j\ge r\}$.

\proclaim{Lemma 4.2} The curves $C^{\pm}(c)$ lie in the regions
$$\left\{z:|\arg(r\inv+z)|\le\pi{cn\ov\al n(r)+cn}\right\}$$
for all $r$ and in $|z+r_2\inv|\ge\dl n\inv$ if $\dl$ is small
enough.
\endproclaim

\demo{Proof} For a point $z$ on either of the curves, say in the upper
half-plane, we have
$$
\align 
c\pi&={\al\ov n}\sum_{j=2}^n\arg(r_j\inv+z)+\arg(z-1)+(c-1)\,\arg z\\
&\ge {\al n(r)\ov n}\arg(r\inv+z)+c\,\arg(r\inv+z),
\endalign 
$$
which gives the first statement of the
lemma. For the second, observe that if $\z =O(\ph_n)$ then
$\si_1'(r_2\inv+\z,c)=\al/n\z+O(1)$. This shows, first, that $u_c^-$ lies to the
right of the circle $|\z|=\dl\,n\inv$ if $\dl$ is small enough and, second, that
$1/\si_1'(z,c)$, thought of a vector, points outward from this circle if $\dl$ is small enough. Since a point of $C^-(c)$ 
moves in the direction of $1/\si_1'(z,c)$ as it moves
away from $u_c^-$ (see (3.7) of [GTW2]), the curve can never pass inside the
circle. Therefore the entire disc $|\z|\le\dl\,n\inv$ lies to the left of
$C^-(c)$. This gives the second statement for $C^-(c)$ and it follows also for
$C^+(c)$ since this is to the right of $C^-(c)$.
$\square$\enddemo

The next lemma, together with (\CH) and the length estimate (\DH), will imply that for $z$ 
large the curves will move in the direction
of $z$ and are well-behaved. If we take any $\bar r<b/(1-b)$ then a
positive proportion of the $r_j$ are greater than $\bar r$ and so by Lemma~4.2 
the
curves lie in a region $$\left\{z:|\arg(\bar r\inv+z)|\le\pi
(1-\dl)\right\}\tag{\regionone}$$ for some $\dl>0$.

\proclaim{Lemma 4.3} We have $z\,\si_1'(z,c)\ra c+\al$ as $n\ra\iy$ and $z\ra\iy$
through region ({\regionone}).
\endproclaim

\demo{Proof} We have $$z\,\si_1'(z,c)=c+\al+O(n\inv)+O(z\inv)+{\al\ov
n}\sum_{j=2}^n{1\ov1+r_jz},$$ and it suffices to show that the last term tends
to 0 as $n\ra\iy$ and $z\ra\iy$ through region ({\regionone}). If $z$ is in this 
region
and $r<\bar r/2$ then $|1+rz|\ge \dl(1+ r|z|)$ for another $\dl$. The same bound
will hold for all $r\le b/(1-b)$ if $z$ is large enough. Choose $M$ large and
break the sum on the right, with its factor $n\inv$, into two parts, the terms
where $r_j|z|<M$ and the terms where $r_j|z|\ge M$. We find that its absolute
value is at most $$n\inv(n-n(M/|z|))+{1\ov \dl M}.$$ The first term tends to 0
as $z\ra\iy$ while the second could have been arbitrarily small to begin with.
$\square$
\enddemo

\noi{\bf Remark}. If $\Pr(p=0)$ is positive then the above has to be modified.
We replace $c+\al$ by $c+\al\,\Pr(p>0)$.

Because of the above lemma we need only consider $z$ in a bounded set. We use the 
fact
that by Lemma~4.2 with $r=r_2$ our curves lie a region
$$\left\{z:|\arg(r_2\inv+z)|\le\pi(1-\dl n\inv),\ \ |r_2\inv+z|\ge\dl
n\inv\right\}. \tag{\regiontwo}$$

\proclaim{Lemma 4.4} For all $z$ in any bounded subset of the region
({\regiontwo}) we have $$|\si_1'(z,c)|\ge \dl\,
n^{-6}\,\left|{(z-u_c^-)\,(z-u_c^+)\ov z(z-1)}\right|$$ for some $\dl>0$ independent of $c$.
\endproclaim
\demo{Proof} To obtain the lower bound we write
$$\Ph(s;z)=\Ph(s_2,\,s_3,\cd,s_n;\,z)={\al\ov n}\sum_{j=2}^n{1\ov s_j+z}+{1\ov
z-1}+ {c-1\ov z}.$$ 
Of course $\si_1'(z,c)=\Ph(r_2\inv,r_3\inv,\cd,r_n\inv)$.
Think of $s_2=r_2\inv$ and $z$ as fixed, and consider the problem of finding
$\inf\,|\Ph(s;z)|$ where $s_3,\cd,s_n$ are subject to the conditions
$$s_j\ge s_2,\ \ \Ph(s;u_c^{\pm})=0.$$
If we take sequences so that the inf is approached
in the limit, then some $s_j$ may tend to infinity, others may tend to $s_2$,
and the rest, if any, tend to values strictly greater than $s_2$. Thus our inf
is equal to the minimum of $|\Ph(s;z)|$, where $\Ph$ now has the form
$$\Ph(s_2,\,s_3,\cd,s_{n'};\,z)={\al\ov n}\sum_{j=2}^{n'}{n_j\ov
s_j+z}+{1\ov z-1}+ {c-1\ov z}$$ with $n'\le n,\ \sum n_j=n-1$, and the $s_j$
with $j>2$ satisfying $s_j>s_2$ and the constraints $\Ph(s;u_c^{\pm})=0.$

Notice that the minimum cannot be zero since $\Ph(s;\,z)$, thought of for the
moment as a function of $z$, has $n'$ finite zeros. It has zeros at $u_c^{\pm}$
and one between each pair of consecutive $-s_j$ since all the coefficients of
$1/(s_j+z)$ are positive. This accounts for all $n'$ zeros, so our $z$ cannot be
one of them.

We apply Lagrange multipliers to find the minimum of $|\Ph(s;z)|^2$ over
$s_3,\cd,s_{n'}$, achieved at
interior points. There are two constraints, hence two multipliers $\lambda$ and
$\mu$. If $p+iq$ is the value $\Ph(s;z)$ where its absolute value achieves its
minimum, then the equations we get are
$$\Re\,(p-iq)\,{1\ov(s_j+z)^2}={\lambda\ov (s_j+u_c^-)^2}+{\mu\ov
(s_j+u_c^+)^2},$$ where we have divided by the factor $n_j$ appearing in all
terms. This is the same sixth degree polynomial equation for all the $s_j$. It
follows that there are at most six different $s_j$. Assuming there are exactly 
six
(if there are fewer the argument is the same and the final estimate is better)
we change notation again and write these as $s_3,\cd,s_{8}$ so that the minimum 
is
achieved for
$$\Ph(s_2,\,s_3,\cd,s_8;\,z)={\al\ov n}\sum_{j=2}^{8}{n_j\ov
s_j+z}+{1\ov z-1}+{c-1\ov z}$$
with other $n_j$.

This has eight zeros. Two of them are $u_c^{\pm}$ and the other six, lying 
between
consecutive $-s_j$, we denote by $u_1,\cd,u_6$. We have the factorization
$$\Ph(s;\,z)={1-c\ov u_c^-\,u_c^+}{(z-u_c^-)\,(z-u_c^+)\ov z(z-1)}
{\prod_{i=1}^6 (1-z/u_i)\ov \prod_{j=2}^{8}(1-z/s_j)},$$ and it remains to find
a lower bound for this. Near $z=0$ we have $\si_1'(z,c)= (1-c)z\inv-1+\al\lan
r\ran +o(1)$, so if $c$ is close to 1 then $(1-c)/u_c^+=1-\al\lan r\ran +o(1)$.
In particular this is bounded away from zero. Thus the first factor above is
bounded away from zero. As for the factors in the products, observe first that
each factor $1-z/s_j$ is bounded since $z$ and all factors $1/s_j$ are. For the
others, we use again the fact that the curves lie in a region ({\regiontwo}). In
any bounded subset of this region each $|1-z/u_i|\ge \eta n\inv$ for some
$\eta>0$. (If $z$ is in a neighborhood of 0 this is clear since each $u_i<0$.
Otherwise write $1-z/u_i=z(z\inv-u_i\inv)$.) Therefore the product of these is
bounded below by a constant times $n^{-6}$. This completes the proof.
$\square$
\enddemo

\vskip-0.2cm Now we can show that the curves $C^{\pm}(c)$ are not too badly 
behaved.

\proclaim{Lemma 4.5} For some constant $A>0$ the length of $C^+(c)$ is 
$O(n^{A})$ 
and
$$\int_{C^-(c)}|z|^{-2}\,|dz|=O(n^{A}).$$
\endproclaim

 \demo{Proof} It follows from Lemma 4.3 that $C^+(c)$ lies in a bounded set. 
For, this
lemma implies that the 
vectors $1/\si_1'(z,c)$ point outward from a large circle 
$|z|=R$, and since by (\CH) $C^+(c)$ goes in the direction opposite to $1/\si_1'(z,c)$,
a point of the curve starting at $u_c^+$ can never pass outside the
circle. Also, some disc $|z|\le\dl(1-c)$ is disjoint from $C^+(c)$ because
$1/\si_1'(z,c)$ points outward from a small enough circle $|z|=\dl(1-c)$ and so
$C^+(c)$ cannot cross into it. It follows that $\si_1'(z,c)$, and so also
$\si_1(z,c)$, is bounded on any portion of $C^+(c)$ close to $z=0$. A similar
argument shows that some disc $|z-1|\le\dl$ lies entirely inside $C^+(c)$.
Finally, we know that $u_c^-$ is within $O(n^{-1/2})$ of $-r_2\inv$ and if
$\z =o(q_1)$ then $\si_1'(r_2\inv+\z,c)=\al/n\z+O(1)$. In particular $u_c^-$
lies in a region $|\z|\ge\dl n\inv$ for some $\dl>0$. Since also 
$\si_1''=-\al/n\z^2+O(1)$,
by Lemma~3.6,
we deduce that $\si_1''(z,c)=O(n)$ when $|z-u_c^-|\le \dl n\inv/2$, thus for 
such $z$
we have $\si_1(z,c)=\si_1(u_c^-,c)+O(n|z-u_c^-|^2)$. But it follows from 
Lemma~4.1
that $\si_1(u_c^-,c)-\si_1(u_c^+,c)>\ph_n^2$, and then, 
since $n\inv=o(\ph_n^2)$, 
$\si_1(u_c^+,c)<\si_1(z,c)$ for $|z-u_c^-|\le \dl n\inv/2$. 
As the maximum of $\si_1(z,c)$ on $C^+(c)$ occurs at $u_c^+$, 
this shows that the distance from $C^+(c)$ to $u_c^-$is at least $\dl n\inv/2$.
With these facts established we use the lower bound of Lemma 4.4,  
the length estimate (\DH) (extended as in the remark following it), and
the obvious upper bound for $|\si_1(z,c)|$ in the region (\regiontwo) to deduce
that the length of $C^+(c)$ is $O(n^{A})$ for some constant $A$.

As for the integral over $C^-(c)$, we observe that, 
since $c<1$ and $cm$ is
an integer, $1-c$ is at least a constant times $n\inv$. Since $C^-(c)$ lies
outside a disc $|z|\le\dl(1-c)$, we have $z\inv=O(n)$ on $C^-(c)$. 
A lower bound
for the distance
from $C^-(c)$ to $u_c^+$ is obtained using the fact that
$\si_1(u_c^-,c)-\si_1(u_c^+,c)>\ph_n^2$.
Since $\si_1'$ is bounded in a neighborhood of $u_c^+$, we have
$\si_1(u_c^-,c)>\si_1(z,c)$ for $|z-u_c^+|$ less than $\ph_n^2$ times a
sufficiently small constant. This shows that $C^-(c)$ is at least this far from
$u_c^+$. We apply the other bounds as before; we think of 
the integral over the portion of
$C^-(c)$ outside a large circle as the sum of 
integrals over the arcs from $a_k$ to $a_{k+1}$ where $a_k$ is
the point of $C^-(c)$ where $|z|=k$. Lemma 4.3 and (\DH) are 
used again here.
$\square$
\enddemo

\vskip-0.2cm
\flushpar{\bf 5. Asymptotic evaluation of the integrals.}
\vskip-0.1cm

We evaluate $I^+(c)$ first when $c-c_n=O(\ph_n)$. Then
$\si_{1zz}(u_{c}^+,c)=\al/\b^2+o(1)$ and so if we set $z=u_{c}^++\z$ we have
$$\si_1(z,c)=\si_1(u_{c}^+,c)+{\al\ov2\b^2}(1+o(1))\z^2$$
as long as $\z=O(\ph_n)$. If $|\z|=\ph_n$ then the real part of the second
term above is less than a negative constant times $\ph_n^2$ and, 
since this real part
decreases as we go out $C^+(c)$, it is at least this negative whenever  
$|\z|\ge 
\ph_n$.
If we recall that this
gets multiplied by $m$ in the exponent and 
the fact that $C^+(c)$ has the length 
at most a
power of $n$ (by Lemma~4.5), we see that the contribution of this part of the 
integral is
$O\left(e^{m\si(u_{c}^+,c)-n\ph_n^2+O(\log n)}\right)$. It follows from Lemma~3.3
and assumption (c)
that with high probability $q_1\gg \log n/n^{1/2}$, and we could have chosen
$\ph_n$ to satisfy this also. Thus, with error
$o(e^{m\si(u_{c}^+,c)})$ the integral $I^+(c)$ is equal to
$${1\ov2\pi 
i}\int_{|\z|<\ph_n}(1+r_1(u_{c}^++\z))\,e^{(n/2\b^2)(1+o(1))\z^2}\,dz\,
e^{m\si_1(u_{c}^+,c)}$$
(since $\al m=n$). Since $\ph_n\gg n^{-1/2}$, in the limit after making the 
variable
change $\z\ra n^{-1/2}\z$
the integration can be taken over $(-i\iy,i\iy)$ (downward really, but we can
reverse the directions of integrations), the linear factor $\z$ contributes 
zero,
and by ({\uest})
$$1+r_1u_c^+=r_1\left(2\b 
n^{-1/2}+{r_1\b^2\ov\al}(c-c_n)+o(n^{-1/2}+|c-c_n|)\right).$$
Thus the integral is asymptotically equal to
$\b\sqrt{2\pi}in^{-1/2}$ times the above and, by ({\sioneest}),
$$
\align
&I^+(c)={r_1\b^2\ov\sqrt{2\pi}}n\inv \left(2 +{r_1\b\ov\al}n^{1/2}
(c-c_n)+o(1+n^{1/2}|c-c_n|\right)\\
&\qquad\qquad
\times \ps_1(-r_1,c)\inv\,e^{-{r_1 \b^2\ov2\al}m\left(c-c_n+
{2\al\ov r_1\b}(1+o(1))n^{-1/2}\right)^2}.
\endalign
$$

This assumed that $c-c_n=O(\ph_n)$. For all $c\ge c_n$ we
use the second part of Lemma~4.1 and again the fact that $C^+(c)$ 
has the length at most a power of $n$. We deduce
$$I^+(c)=O\left(\ps_1(-r_1\inv,c)\,e^{-\eta n^{1/2}\,(c-c_n)+O(\log n)}\right)$$
for $c\ge c_n$.

For the integral over $C^-(c)$ we use the last part of Lemma~4.1 and the
second part of Lemma~4.5. These imply that the integral over $C^-$ is
$$O\left(\ps_1(-r_1\inv,c)\inv\,e^{-n\ph_n^2+O(\log 
n)}\right)=o(\ps_1(-r_1\inv,c)).$$

But our integral for $I^-(c)$ is {\it not} taken over $C^-(c)$. Recall that the 
original
contour must have all the $-r_j\inv$
on the outside whereas $-r_1\inv$ is inside (more precisely, on the other side 
of) $C^-(c)$.
Therefore if we deform the contour to
$C^-(c)$ we pass through the pole at $-r_1\inv$. Thus
$$I^-(c)=r_1\,\psi_1(-r_1\inv,c)\inv+o(\ps_1(-r_1\inv,c)).$$

Now recall that in $I^+(c)$ we set $c-c_n=h'+j+\ell$, in $I^-(c)$ we set 
$c-c_n=h'+\ell+k$
and then we sum over $\ell$ to get the matrix product. Recall also that
$\ps_1(-r_1\inv,c)=\ps_1(-r_1\inv)\,(-r_1)^{-m(c-c_n)}$. The
factors $(-r_1)^{-m(c-c_n)}$ in $I^+(c)$ and $(-r_1)^{m(c-c_n)}$ in $I^-(c)$ will
combine to give $(-r_1)^{m(k-j)}$ which can be eliminated without affecting the
determinant. It follows that we can modify the expressions for $I^\pm(c)$ by 
removing
these factors. We can also remove the factors $\ps_1(-r_1\inv)^{\pm1}$ since they 
cancel
upon multiplying. Thus our replacements are
$$I^+(c)\ra {r_1\b^2\ov\sqrt{2\pi}}n\inv \left(2 +{r_1\b\ov\al}n^{1/2}
(c-c_n)\right)e^{-{r_1 \b^2\ov2\al}m\left(c-c_n+
{2\al\ov r_1\b}(1+o(1))n^{-1/2}\right)^2},$$
if $c-c_n=O(\ph_n)$, 
and 
$$I^+(c)\ra O\left(e^{-\eta n^{1/2}(c-c_n)+O(\log n)}\right),$$
if $c>c_n.$
Furthermore, $I^-(c)\ra r_1+o(1).$

Recall next that we set $h'=sn^{1/2}$ and in $I^+(c),\ 
c=c_n+sn^{1/2}+\lfloor xn^{1/2}\rfloor+\lfloor zn^{1/2}\rfloor$,
so that
$$c-c_n=(s+x+z+o(1))n^{1/2}/m=\al(s+x+z+o(1))n^{-1/2},$$
and eventually we multiply by $n$ because of the scaling.
Take first the case $c-c_n=O(\ph_n)$, that is, $x+z=O(n^{1/2}\ph_n)$.
Since $m=n/\al$ and $r_1\,\b=\tau\inv\,(1+o(1))$ the modified $I^+(c)$ equals
$${r_1^2\b^3\ov\sqrt{2\pi}}n\inv (2\tau +s+x+z+o(1+x+y))
e^{-(2\tau +s+x+z+o(1))^2/2\tau^2}.$$
On the other hand, $I^-(c)$ is equal to $r_1$ with error $o(1)$.
The result of multiplying these together, multiplying by $n$, and integrating
with respect to $z$ over $(0,\,\infty)$, is asymptotically equal to
$${1\ov\sqrt{2\pi}\tau}e^{-(2\tau +s+x)^2/2\tau^2}.\tag{\kernel}$$
This holds for $c-c_n=O(\ph_n)$. If $c-c_n\ge \ph_n$ 
we have, for our modified $I^+(c)$, the estimate
$$O\left(e^{-\eta n^{1/2}(c-c_n)+O(\log n)}\right)=O(n\inv).$$
Integrating the square of this over a region $x+z=O(n^{1/2})$ will give $o(1)$.

It follows that the matrix product scales to the operator on $(0,\,\iy)$ with 
kernel
({\kernel}). This is a rank one kernel so its Fredholm determinant equals one
minus its trace, which equals
$${1\ov\sqrt{2\pi}\tau}\int_{-\iy}^{2\tau+s}e^{-x^2/2\tau^2}.$$
This establishes the convergence in probability statement of Theorem 1.\sp

\vskip-0.1cm
\noi{\bf Remark}. One could rightly object that to scale a product to a trace 
class operator
we should know
that each factor scales in Hilbert-Schmidt norm. In our case the second limiting 
kernel
is a constant and the product is not even Hilbert-Schmidt. But we could have
multiplied the kernel of the first operator by $(1+x)\,(1+z)$ and the kernel
of the second operator by $(1+z)\inv\,(1+y)\inv$. This would not have affected 
the
determinant of the product, both operators would have scaled in Hilbert-Schmidt 
norm
and the product would have scaled in trace norm to the rank one kernel
$${1\ov\sqrt{2\pi}\tau}\,e^{-(2\tau +s+x)^2/2\tau^2}{1+x\ov 1+y}$$
which has the same Fredholm determinant.\sp
 
\vfill\eject 

\subheading{6. Almost sure convergence}

What is needed, and all that is needed, is an ``almost sure'' substitute for
Lemma 3.6 under assumptions (a$'$) and (b$'$). We begin with a lemma on extreme order 
statistics of uniform random variables, 
part or all of which may well be in the literature.
\flushpar{\bf Lemma 6.1} {\it Let $a>1$
be arbitrary. Then, almost surely, 
$$t_1\ge {\eta\ov n\,\log^an},\ \ \ \ {t_1\ov t_2}\le 1-{1\ov \log^{a}n},$$
for sufficiently large $n$. Here, 
$\eta$ is a positive constant depending on $a$.}
 
\comment
\vskip-0.5cm\proclaim{Lemma 6.1} Let $a>1$
be arbitrary. Then, almost surely, 
$$t_1\ge {\eta\ov n\,\log^an},\ \ \ \ {t_1\ov t_2}\le 1-{1\ov \log^{a}n},$$
for sufficiently large $n$. Here, 
$\eta$ is a positive constant depending on $a$.
\endproclaim
\endcomment

\demo{Proof} We use the notation $t_{n,j}$ for our $t_j$ to display
their dependence on $n$. We have
$$\Pr(t_{n,1}\le\dl)=1-(1-\dl)^n\sim n\dl\ \ {\text {if}}\ n\dl=o(1).$$
In particular
$$\Pr\left(t_{2^k,1}\le {2^{-k}\ov k^a}\right)\sim {1\ov k^a}.$$
It follows that, a.s. for sufficiently large $k$ we have
$$t_{2^k,1}>{2^{-k}\ov k^a}.$$
Take any $n$ and let $k$ be such that $2^{k-1}<n\le 2^k$. From the above 
we have, 
a.s.
for sufficiently large $n$
$$t_{n,1}\ge t_{2^k,1}>{2^{-k}\ov k^a}\ge{\eta\ov n \log^an},$$
for some $\eta$.

For the ratio we use the fact that 
$$\Pr\left({t_{n,j}\ov t_{n,j+1}}>1-\dl\right)=
1-(1-\d)^j\sim j\dl\ \ {\text {if}}\ 
j\dl=o(1).\tag{6.1}$$
Now suppose that 
$${t_{n,1}\ov t_{n.2}}> 1-{1\ov \log^{a}n}\tag{6.2}$$
and let $k$ be such that $2^{k-1}<n\le 2^k$. 
Take any $J$ (which will eventually be of 
order $\log k$). Then there are two possibilities:

\beginitems

\item{(1)} $t_{2^k,j}\le t_{n,1}$ for all $j\le J$;

\item{(2)} $t_{2^k,j}>t_{n,1}$ for some $j\le J$.

\enditems

Consider possibility (1) first. Let $G_n$ be 
the event that $\ t_{n,1}\le a\log\log n/n$. 
By Ex.~4.3.2 of [Gal], $P(G_n$ eventually$)=1.$
Moreover, 
$$
\align 
&\Pr(\{t_{2^k,j}\le t_{n,1}\ {\text {for\ all}}\ j\le J\}\cap G_n)
\le \Pr(t_{2^k,j}\le 2\log\log n/n\ {\text {for\ all}}\ j\le J)\\
&\le{{2^k}\choose{J!}}\left(2\,{\log\log n\ov n}\right)^J
\le e^{J\log\log k-J\log J+AJ}, 
\endalign $$
for some constant $A$. If $J=B\log k$ then the bound above 
equals $e^{-B(\log B-A)\log k}$,
so if we choose $B$ large enough the sum over $k$ of these probabilities will
be finite. With this $J$, (1) can therefore a.s. 
occur for only finitely many $k$. 

Next consider possibility (2) 
and let $j$ be the smallest integer $\le J$ such that
$t_{2^k,j}>t_{n,1}$. Then $t_{2^k,j}\le t_{n,2}$ and $t_{n,1}=t_{2^k,\ell}$ for 
some $\ell<j$. It follows that 
$t_{2^k,j-1}/t_{2^k,j}> t_{n,1}/t_{n,2}$ and by (6.2) this is 
at least $1-C/k^{a}$, for some constant $C$ (which will 
change from appearance to appearance). 
Therefore, by (6.1), 
$$
\align 
&P((6.2)\text{ and }(2)\text{ both happen})\\
&\le 
P(t_{2^k,j-1}/t_{2^k,j}>1-C/k^a\text{ for some }j\le J)\le 
C J^2/k^a\le C\log^2k/k^a. 
\endalign
$$ 
It follows that (2) and (6.2) can 
happen together only for
finitely many $n$. The upshot is that a.s. the inequality (6.2) can occur for
only finitely many $n$, which completes the proof.
$\square$
\enddemo

We are now ready to prove our substitute for Lemma~3.6. Recall that
we can set $q_j=G^{-1}(t_j)$. The assumption (a$'$) implies that 
$G$ is continuous near 0, so that $G(G^{-1}(x))=x$ 
for small $x$.

\proclaim{Lemma 6.2} Suppose (a$'$) and (b$'$) are satisfied. 
Then there 
exists a sequence $\ph_n\gg \log n/n^{1/2}$ such 
that a.s. for any sequence $\{v_n\}$ lying in the disc with diameter the real 
interval
$[-r_1\inv-O(\ph_n),\,\xi]$ we have 
$$\lim_{n\ra\iy}{1\ov 
n}\sum_{j=2}^n{r_j\ov(1+r_jv_n)^2}=\lan{r\ov(1+r\xi)^2}\ran.$$
\endproclaim

\demo{Proof} From 
the proof of Lemma 3.6 we see that we want to show that, for
some sequence $\ph_n$ as described, we have a.s.
$$\lim_{n\ra\iy}{1\ov n}\sum_{j=2}^n{q_1\ov q_j\,(q_j-(q_1+O(\ph_n))^2}=0.$$
Assumption (a$'$) implies that  
$$
\frac xy\ge \(\frac{G^{-1}(x)}{G^{-1}(y)}\)^\ga,
$$
when $x\le y$ are small enough. Therefore, it follows from 
the second part of Lemma 6.1, that a.s. for large 
$n$, 
$${q_1\ov q_2}\le 1-{\eta\ov \log^{a}n}\tag{6.3}$$
for another constant $\eta>0$. Set 
$$\ps_n={1\ov2}{\eta\ov \log^{a}n}q_2.$$

Let us show that $\ps_n\gg \log n/n^{1/2}$.
Assumption (a$'$) implies that  
$G^{-1}(x)$ is at most a
constant times  $x^{1/\g}$, 
thus the fact that $t_1=O(\log\log n/n)$ shows that
$q_1$ is at most a constant times $(\log\log n/n)^{1/\ga}$.
Furthermore, 
assumption
(b$'$) gives, with a slightly smaller $\nu$,
$x^2\gg G(x)\log^{\nu} x\inv$. 
Applying this with 
$x=q_1=G^{-1}(t_1)$ and using the first part of Lemma 6.1 gives 
$$q_1^2\gg {1\ov  n \log^a n}\log^{\nu}q_1\inv.$$
We therefore deduce that
$$q_1^2\gg {1\ov  n} \log^{\nu-a}n\tag{6.4}$$
for a slightly smaller $\nu$ than in (b$'$).
By (6.3), the same holds for $q_2$ and so 
$$\ps_n^2\gg{1\ov n}\log^{\nu-3a}n$$
and $\ps_n\gg \log n/n^{1/2}$ 
as long as $\nu-3a>2$. Since $a>1$ is arbitrary the 
requirement becomes $\nu>5$. But from (a$'$) and (b$'$) we 
see that necessarily 
$\ga>2$, so that $\nu>8$.

If $j\ge 2$, then (6.3) and the inequality
$q_2\le q_j$  imply that
$$q_j-(q_1+\ps_n)\ge{1\ov2}{\eta\ov \log^{a}n}q_j.$$

We take for $\{\ph_n\}$ any sequence satisfying
$${\log n\ov n^{1/2}}\ll\ph_n\ll\ps_n.$$
At this point we follow the proof of Lemma 3.6 
to see that the expression
$${\log^{2a}n\ov n}\sum_{j=2}^n{q_1\ov q_j^3}\tag{6.5}$$
needs to go to 0 a.s. to conclude the proof of this lemma.
This is what we will demonstrate. 

For any $k_n$, 
if we separate the sum in (6.5) over $j\le k_n$ from the sum over $j>k_n$,
we see that (6.5) is at most
$${\log^{2a}n\ov n\,q_1^2}\,k_n+\log^{2a}n\,{q_1\ov q_{k_n+1}}\,
{1\ov n}\sum_{j=1}^n{1\ov q_j^2}.\tag {6.6}$$

We first determine $k_n$ so the second term in (6.6) goes a.s. to 0. 
By strong law, $n^{-1}\sum q_j\to \<q^{-2}\>$ a.s., so 
$\log^{2a}n\,{q_1/q_{k_n+1}}$ needs to go to 0. 
We have, for each $\dl>0$, 
$$
\align
&\Pr \left(\log^{2a}n{q_1\ov q_{k_n+1}}\ge\dl\right)=
\Pr \left({G^{-1}(t_1)\ov G^{-1}(t_{k_n+1})}\ge{\dl\ov \log^{2a}n}\right)\\
&\le \Pr \left({t_1\ov t_{k_n+1}}\ge\left({\dl\ov \log^{2a}n}\right)^\ga\right)
=\left(1-\left({\dl\ov \log^{2a}n}\right)^\ga\right)^{k_n}\\
&\le e^{-\left({\dl\ov \log^{2a}n}\right)^\ga k_n}.
\endalign
$$
This is summable over $n$ if we choose
$$k_n=\lfloor\log^a n\left({\log^{2a}n}\right)^\ga\rfloor+1.$$
With this choice, the second summand in (6.6) therefore goes to 0 a.s. 

On the other hand, the first term in (6.6) is
with the same choice of $k_n$ at most a constant times
$${\log^{(2\ga +3)a}n\ov n\,q_1^2},$$
and from (6.4) this is $o(1)$ times $\log^{(2\ga+4)a-\nu}n$. Since $a>1$
was arbitrary and $\nu>2\ga+4$, we can make $(2\ga+4)a-\nu<0$ 
and then the first summand in (6.6) goes to 0 a.s. This completes the proof.
$\square$
\enddemo

With this lemma in place of Lemma 3.6 the reader will find 
that all subsequent 
limits and estimates in  Sections 4 and 5
will hold almost surely, thus giving the second 
statement of the 
theorem. The reason our sequence
had to satisfy $\ph_n\gg \log n/n^{1/2}$ is that errors of the form
$O\left(e^{-n\ph_n^2+O(\log n)}\right)$ appeared in the evaluation of 
$I^{\pm}(c)$
and these had to be $o(1)$.

\vfill\eject

\vskip1cm

\centerline{REFERENCES}

\vskip0.5cm
  
\item{\bf [BCKM]}
J.-P.~Bouchaud, L.~F.~Cugliandolo, J.~Kurchan, M.~M\'ezard, {\it 
Out of equilibrium dynamics in spin--glasses and 
other glassy systems.\/}  In ``Spin Glasses and 
Random Fields,'' A.~P.~Young, editor, World Scientific, 1998.

\item{\bf [BDJ]}
J.~Baik, P.~Deift, K.~Johansson, {\it 
On the distribution of the length of the 
longest increasing subsequence
of random permutations.\/} J.\ Amer.\ Math.\ Soc. 12 (1999),
1119--1178.

\comment

\item{\bf [BFL]}  I.~Benjamini, P.~A.~Ferrari, C.~Landim, {\it 
Asymptotic conservative processes with random rates.\/}  
Stochastic Process.\ Appl.\ 61 (1996), 181--204. 

\endcomment
 
\item{\bf [BR]} 
J.~Baik, E.~M.~Rains, {\it Limiting distributions for a polynuclear
growth model with external sources.\/}
J.\ Statist.\ Phys.\ 100 (2000), 523--541.
 
\item{\bf [FIN1]} L.~R.~G.~Fontes, M.~Isopi, C.~M.~Newman, {\it 
Random walks with strongly inhomogeneous rates 
and singular diffusions: convergence, 
localization and aging in one dimension.\/} 
Preprint (ArXiv: math.PR/0009098). 

\item{\bf [FIN2]} L.~R.~G.~Fontes, M.~Isopi, C.~M.~Newman, {\it 
Chaotic time dependence in a disordered spin system.\/}  
Probab.\ Theory Relat.\ Fields 115 (1999), 417--443.

\item{\bf [FINS]} L.~R.~G.~Fontes, M.~Isopi, C.~M.~Newman, 
D.~L.~Stein, 
{\it Aging in 1D Discrete Spin Models and Equivalent Systems.\/}  
Phys.\ Rev.\ Lett.\ 87 (2001). 110201--1.

\item{\bf [Gal]} J.~Galambos, ``The Asymptotic Theory of Order Statistics.''
Second edition. Krieger, 1987.

\item{\bf [GG]} J.~Gravner, D.~Griffeath, {\it Cellular automaton growth on 
$Z^2$: theorems, examples, and
problems.\/} Adv.\ in Appl.\ Math.\ 21 (1998), 241--304.

\item{\bf [Gri1]} D.~Griffeath,  ``Additive and Cancellative 
Particle Systems.'' Lecture Notes in Mathematics
724, Springer, 1979.

\item{\bf [Gri2]} D.~Griffeath, {\it Primordial 
Soup Kitchen.\/} {\tt psoup.math.wisc.edu}

\item{\bf [Gra]} J.~Gravner, {\it Recurrent ring 
dynamics in two--dimensional excitable cellular 
automata.\/} J.\ Appl.\ Prob.\ 36 (1999), 492--511. 

\item{\bf [GTW1]} J.~Gravner, C.~A.~Tracy, H.~Widom, 
{\it Limit theorems for height fluctuations in a class of discrete space
and time growth models. \/} 
J.\ Statist.\ Phys.\ 102 (2001), 1085--1132.

\item{\bf [GTW2]} J.~Gravner, C.~A.~Tracy, H.~Widom, 
{\it  A growth model in a random environment.\/} To 
appear in Ann.\ Probab.\ (ArXiv: math.PR/0011150).

\item{\bf [Joh1]}
K.~Johansson, {\it Shape fluctuations and random matrices.\/} 
Commun.\ Math.\ Phys.\ 209 (2000), 437--476.

\item{\bf [Joh2]} K.~Johansson, {\it 
Discrete orthogonal polynomial ensembles and the Plancherel
measure.\/} Ann.\ Math.\ 153 (2001), 259--296. 

\item{\bf [Lig]} T.~Liggett, ``Interacting Particle Systems.'' 
Springer--Verlag, 1985. 

\item{\bf [Mea]} P.~Meakin, ``Fractals, scaling and growth 
far from equilibrium.'' 
Cambridge University Press, 1998.

\item{\bf [MPV]} M.~M\'ezard, G.~Parisi, M.~A.~Virasoro, 
``Spin Glass Theory and Beyond.'' 
World Scientific, 1987.

\item{\bf [NSt1]} C.~M.~Newman, D.~L.~Stein,   
{\it Equilibrium pure states and nonequilibrium chaos.\/} 
J.\ Statist.\ Phys.\ 94 (1999), 709--722.

\item{\bf [NSt2]} C.~M.~Newman, D.~L.~Stein, {\it 
Realistic spin glasses below eight dimensions: a highly disordered
view.\/} Phys. Rev. E (3) 63 (2001), no. 1, part 2, 016101, 9 pp. 

\item{\bf [NSv]} P.~Norblad, P.~Svendlindh, 
{\it Experiments on spin glasses.\/} In ``Spin Glasses and 
Random Fields,'' A.~P.~Young, editor, World Scientific, 1998. 

\comment

\item{\bf [NV]} C.~M.~Newman, S.~B.~Volchan, {\it 
Persistent survival of one-dimensional contact 
processes in random environments. \/} 
Ann.\ Probab.\ 24 (1996), 411--421.

\endcomment

\item{\bf [PS]} M.~Pr\"ahofer, H.~Spohn, {\it Universal distribution 
for growth processes in $1+1$ dimensions and random 
matrices\/}. Phys.\ Rev.\ Lett.\ 84 (2000), 
4882--4885.

\comment

\item{\bf [PS2]} 
M.~Pr\"ahofer, H.~Spohn, {\it Scale Invariance of the PNG Droplet
and the Airy Process.\/}  Preprint (ArXiv: math.PR/0105240).

\item{\bf [Rai]} E.~M.~Rains, {\it A mean identity 
for longest increasing subsequence problems.\/}
Pre\-print (arXiv: math.CO/0004082).

\item{\bf [Sep1]} T.~Sepp\" al\" ainen, {\it Increasing sequences of 
independent points on the planar lattice.\/} 
Ann.\ Appl.\ Probab.\ 7 (1997), 886--898.
 
\item{\bf [Sep2]} T.~Sepp\"al\"ainen, {\it 
Exact limiting shape for a simplified model of first-passage
percolation on the plane.\/} Ann.\ Probab.\ 26 (1998), 1232--1250.

\endcomment

\item{\bf [SK]} T.~Sepp\"al\"ainen,  J.~Krug, {\it Hydrodynamics and 
platoon formation for a  totally asymmetric
exclusion model with particlewise disorder.\/} 
J.\ Statist.\ Phys.\ 95 (1999), 525--567. 

\item{\bf [Sos]} A.~Soshnikov, {\it  
Universality at the edge of the spectrum in Wigner random matrices.\/} 
Commun.\ Math.\ Phys.\ 207 (1999), 697--733.

\item{\bf [Tal]} M.~Talagrand, 
{\it Huge random structures and mean field models for spin glasses.\/}
Doc.\ Math., Extra Vol. I (1998), 507--536. 

\item{\bf [TW1]} C.~A.~Tracy, H.~Widom, 
{\it Level spacing distributions and the Airy kernel.\/}
Commun.\ Math.\ Phys.\ 159 (1994), 151--174.

\item{\bf [TW2]} C.~A.~Tracy, H.~Widom,
{\it Universality of the Distribution Functions of 
Random Matrix Theory. II.\/} In 
``Integrable Systems: From Classical to Quantum,'' 
J.~Harnad, G.~Sabidussi and P.~Winternitz, editors, 
American Mathematical Society, Providence, 2000. Pages 251--264. 
 
\enddocument